\documentclass[a4paper,11pt]{amsart}

\usepackage{geometry}
\usepackage{amssymb,mathrsfs,amsmath,amsthm}
\usepackage{enumitem}
\usepackage{autobreak}
\allowdisplaybreaks
\usepackage{bm}
\usepackage{bbm}

\usepackage{graphicx}
\usepackage{color}
\usepackage{booktabs}
\usepackage{hyperref}
\usepackage{cleveref}

\newtheorem{theorem}{Theorem}[section]
\newtheorem{lemma}{Lemma}[section]
\newtheorem{proposition}{Proposition}[section]
\newtheorem{corollary}{Corollary}[section]

\newtheorem{remark}{Remark}[section]

\theoremstyle{definition}
\newtheorem{example}{Example}[section]

\numberwithin{equation}{section}
\numberwithin{figure}{section}
\DeclareMathOperator{\TV}{TV}
\DeclareMathOperator{\Per}{Per}
\DeclareMathOperator{\Lip}{Lip}
\DeclareMathOperator{\diver}{div}
\DeclareMathOperator{\id}{id}
\DeclareMathOperator*{\infess}{ess\,inf}
\newcommand{\R}{\mathbb R}
\newcommand{\A}{\mathcal A}
\newcommand{\K}{\mathcal K}
\newcommand{\Prob}{\mathcal P}
\newcommand{\dd}{\,\mathrm d}
\newcommand{\ind}{\mathbf 1}
\numberwithin{equation}{section}
\makeatletter
\@namedef{subjclassname@2020}{\textup{2020} Mathematics Subject Classification}
\makeatother

\begin{document}

\title{Total Variation Regularized Barycenters in Wasserstein Space}

\author{Qing Zhao}
\address{School of Mathematics and Shing-Tung Yau Center, Southeast University, Nanjing 211189, P. R. China}
\email{zhao-qing@seu.edu.cn}

\author{Ziwei Nie$^*$}
\address{School of Mathematics, Nanjing University, Nanjing 210093, China}
\email{nieziwei@nju.edu.cn}


\date{\today}
\thanks{$^*$ Corresponding author: Ziwei Nie.}



\begin{abstract}
We study total variation regularization of Wasserstein
barycenters on bounded connected Lipschitz domains.
We prove existence for arbitrary prescribed probability
measures and uniqueness when at least one is absolutely
continuous. Our main results are attained dual formulations
in terms of relative perimeters and vector fields with
zero normal trace. These yield necessary and sufficient
optimality conditions and a variational characterization
of superlevel sets. We also establish a maximum principle
on convex domains, stability of minimizers, and an exact
reduction to a single TV-JKO step in one dimension.
The existence, uniqueness, and perimeter duality results
extend to costs $|x-y|^p$, $p>1$, under the same assumptions.
\end{abstract}

\maketitle

\section{Introduction}
Wasserstein barycenters provide a notion of averaging
probability measures that accounts for the geometry
of their underlying space.
Agueh and Carlier~\cite{bib1} developed their variational
theory and established connections with multi-marginal
optimal transport.
Subsequent work has extended this framework to more
general metric settings and studied its applications
in imaging and statistics; see, for example,
\cite{bib8,bib9,bib10,bib16,bib21}. Regularization provides a way to impose additional
structure on the barycenter.
Bigot, Cazelles and Papadakis~\cite{bib6} studied
penalized Wasserstein barycenters, while
Carlier, Eichinger and Kroshnin~\cite{bib19}
investigated entropy penalties on the barycenter density.
Such penalties act on the unknown measure itself.
A different approach regularizes the transport problems,
for example by adding entropy to the transport plans;
see~\cite{bib22}. Chizat~\cite{Chizat26} studies a doubly regularized
formulation combining entropic transport costs with
an entropy penalty on the barycenter, and establishes
regularity and stability results.
Our formulation retains the unregularized transport
costs and penalizes the total variation of the
barycenter density.

Let $\Omega\subset\R^d$ be a nonempty bounded connected
Lipschitz domain, and let
$\mu_1,\ldots,\mu_n\in\Prob(\overline\Omega)$. We seek a barycenter whose density has bounded variation
and penalize its spatial variation.
To formulate this problem on $\Prob(\overline\Omega)$,
we extend the total variation of densities by
$+\infty$ outside this class. For $\nu\in\Prob(\overline\Omega)$, define
$$
\mathcal T(\nu):=
\begin{cases}
|D\rho|(\Omega),
& \text{if }\nu=\rho\,dx
  \text{ for some }\rho\in BV(\Omega),\\
+\infty,
& \text{otherwise}.
\end{cases}
$$
We consider
$$
\inf_{\nu\in\Prob(\overline\Omega)}
\left\{
\sum_{i=1}^n W_2^2(\mu_i,\nu)+\mathcal T(\nu)
\right\}.
$$
The probability measure with constant density
$1/|\Omega|$ on $\Omega$ has finite energy,
so the infimum is finite. By the definition of $\mathcal T$, the problem
is therefore equivalent to minimization over
probability densities in $BV(\Omega)$:
\begin{equation}\label{eq:intro-problem}
\inf_{\rho\in\A}
\left\{
\sum_{i=1}^n W_2^2(\mu_i,\rho\,dx)
+|D\rho|(\Omega)
\right\},
\qquad
\A=
\left\{
\rho\in BV(\Omega):
\rho\geq0\ \text{a.e.},\
\int_\Omega\rho\,dx=1
\right\}.
\end{equation}
The prescribed measures may be singular, whereas
every minimizer has a density in $BV(\Omega)$.
The total variation penalty controls spatial variation
while allowing jump discontinuities.
Its use in image restoration goes back to
Rudin, Osher and Fatemi~\cite{ROF92}.
In the present setting, it controls the spatial
variation of the barycenter without requiring
its density to be smooth or strictly positive.

The closest related variational problem is the TV-JKO scheme studied by Carlier and Poon~\cite{CP19}.
For a given probability measure $\mu$, one step
of this scheme minimizes
$$
\frac{1}{2\tau}W_2^2(\mu,\rho\,dx)+|D\rho|(\Omega)
$$
over $\rho\in\A$.
Thus \eqref{eq:intro-problem} with $n=1$ corresponds
to $\tau=1/2$. We adapt the entropic approximation of Carlier and
Poon~\cite{CP19} to derive vector-field optimality
conditions for TV-regularized Wasserstein barycenters,
allowing singular prescribed measures.
We also extend their projection argument to obtain
a maximum principle on convex domains. We further establish an attained dual formulation
through relative-perimeter inequalities and relate it
to the vector-field formulation.
This duality characterizes almost every superlevel set
of positive measure of a minimizing density as a
minimizer of a perimeter ratio.
In one dimension, an exact variance identity reduces
the barycenter problem to a single TV-JKO step.

\subsection{Main results}

We first prove existence for arbitrary
$\mu_1,\ldots,\mu_n\in\Prob(\overline\Omega)$.
If at least one of these measures is absolutely
continuous, the classical strict-convexity argument
for the transport cost gives uniqueness.
We then establish an attained dual formulation
involving transport potentials, a scalar multiplier,
and relative-perimeter constraints.

The duality identity yields necessary and sufficient
conditions for optimality.
These conditions characterize optimal transport plans
and imply that almost every superlevel set of positive
measure of a minimizing density minimizes the
relative-perimeter ratio associated with the sum
of the optimal transport potentials.
They also give a transport-map representation
when a prescribed measure is absolutely continuous
and a criterion for the uniform density to be optimal.

Following the entropy approximation in~\cite{CP19},
we establish a second attained dual formulation
in terms of bounded vector fields.
The resulting optimality conditions retain the
multiplier for the mass constraint and a nonnegative
multiplier associated with the constraint $\rho\geq0$.
This construction does not require the prescribed
measures to be absolutely continuous.

We further prove a maximum principle on convex domains
and stability under approximation of the prescribed
measures.
In one dimension, a variance identity reduces the
problem exactly to a TV-JKO step with the ordinary
Wasserstein barycenter as reference measure.
Finally, existence, uniqueness under the same
absolute-continuity condition, and the perimeter
dual formulation extend to costs $|x-y|^p$ with $p>1$.

\subsection{Organization of the paper}

Section~\ref{sec2} introduces the notation
and recalls the analytical tools.
Section~\ref{sec:duality} establishes existence,
uniqueness, and the perimeter dual formulation.
Section~\ref{sec:optimality} derives the corresponding
optimality conditions and transport representations.
Section~\ref{sec5} develops the entropy
approximation and proves vector-field duality.
Section~\ref{sec:qualitative} contains the maximum
principle, the one-dimensional reduction, and
the stability result.
Section~\ref{sec:extensions} treats the extension
to $p$-Wasserstein costs.

\section{Notation and preliminaries}\label{sec2}
We introduce the notation used throughout the paper
and recall the compactness, optimal transport, and minimax results needed below.

\subsection{Notation}
Throughout, $\Omega\subset\R^d$, $d\geq1$, is a nonempty
bounded connected open set with Lipschitz boundary.
Set $
K=\overline\Omega$ and $D=\operatorname{diam}(K)$. Convexity of $\Omega$ will be assumed only when explicitly
stated. We denote by $\mathcal L^d$ the Lebesgue measure, and use $|E|$ and $\ind_E$ for the Lebesgue measure and
characteristic function of a measurable set $E$.
The notation $|\cdot|$ also denotes the Euclidean norm,
with the meaning determined by the argument.
Unless another measure is specified, almost-everywhere
statements refer to Lebesgue measure.

For $1\leq s<\infty$, we use the norms
$$
\|u\|_{L^s(\Omega)}
=\left(\int_\Omega |u|^s\,dx\right)^{1/s},
\qquad
\|u\|_{L^\infty(\Omega)}
=\operatorname*{ess\,sup}_{x\in\Omega}|u(x)|.
$$
The same notation is used for vector-valued functions,
with the Euclidean norm inside the integral.
We write $C(K)$ for the real-valued continuous functions
on $K$, equipped with the uniform norm
$$
\|f\|_{C(K)}=\max_{x\in K}|f(x)|.
$$
The abbreviation $\|\cdot\|_\infty$ denotes the uniform
norm for functions in $C(K)$ and the essential supremum
norm for measurable functions on $\Omega$. We use $C_c^1(\Omega;\R^d)$ for continuously differentiable
vector fields with compact support in $\Omega$. For a real-valued function $u$, set
$$
u_+=\max\{u,0\},
\qquad
u_-=\max\{-u,0\},
\qquad
\overline u=\frac1{|\Omega|}\int_\Omega u\,dx
\quad\text{if }u\in L^1(\Omega).
$$
The symbols $\to$, $\rightharpoonup$, and
$\rightharpoonup^*$ denote strong, weak, and weak-star
convergence in the indicated function spaces. For $u\in L^1(\Omega)$, its distributional gradient is
denoted by $Du$. Define
\begin{equation}\label{eq:TV}
\TV(u)
=
\sup\left\{
\int_\Omega u\,\diver z\,dx:
z\in C_c^1(\Omega;\R^d),\ \|z\|_\infty\leq1
\right\},
\end{equation}
where $
\diver z=\sum_{j=1}^d\partial_j z_j$. Then $
BV(\Omega)=\{u\in L^1(\Omega):\TV(u)<\infty\}$.
For $u\in BV(\Omega)$, $Du$ is a finite vector-valued
Radon measure, and
$$
\TV(u)=|Du|(\Omega),
\qquad
\|u\|_{BV(\Omega)}
=\|u\|_{L^1(\Omega)}+|Du|(\Omega).
$$
Here $|Du|$ denotes the total variation measure of $Du$.
For functions outside $BV(\Omega)$, the expression
$|Du|(\Omega)$ is understood as $+\infty$.
For a measurable set $E\subset\Omega$, its relative
perimeter is $
\Per(E;\Omega)=|D\ind_E|(\Omega)$. Thus only variation inside $\Omega$ is counted.

We write $\Prob(K)$ for the Borel probability measures
on $K$. For $p\geq1$, let
$$
\Prob_p(K)
=
\left\{\mu\in\Prob(K):
\int_K|x|^p\,d\mu(x)<\infty\right\}.
$$
Since $K$ is compact, $\Prob_p(K)=\Prob(K)$.
For $\mu,\nu\in\Prob(K)$, the set $\Pi(\mu,\nu)$ consists
of probability measures on $K\times K$ with first marginal
$\mu$ and second marginal $\nu$.
For a Borel map $T$, its pushforward is defined by $
(T_\#\mu)(B)=\mu(T^{-1}(B))$
for every Borel set $B$ in the target space.
We denote the identity map by $\id$ and the Dirac
measure at $x$ by $\delta_x$.
Narrow convergence $\mu_k\rightharpoonup\mu$ means
$$
\int_K\varphi\,d\mu_k\rightarrow\int_K\varphi\,d\mu
\qquad\text{for every }\varphi\in C(K).
$$

\subsection{Optimal transport and analytical preliminaries}

We first recall the compactness and coarea properties
of $BV$ functions, followed by the transport duality
and minimax results used in the analysis.

On a bounded Lipschitz domain, every bounded sequence
in $BV(\Omega)$ has a subsequence converging strongly
in $L^1(\Omega)$ to a function in $BV(\Omega)$.
Moreover,
$$
u_k\to u\text{ in }L^1(\Omega)
\quad\Longrightarrow\quad
|Du|(\Omega)\leq\liminf_{k\to\infty}|Du_k|(\Omega).
$$
The embedding $BV(\Omega)\hookrightarrow L^s(\Omega)$
is compact for $1\leq s<d/(d-1)$ when $d\geq2$,
and for every finite $s\geq1$ when $d=1$.
We also use the composition estimate $|D(\Phi\circ u)|(\Omega)
\leq\Lip(\Phi)|Du|(\Omega)$
for $u\in BV(\Omega)$ and Lipschitz functions
$\Phi:\R\to\R$. For nonnegative $u\in BV(\Omega)$, the coarea
and layer-cake formulas give
\begin{align}
|Du|(\Omega)
&=\int_0^\infty\Per(\{u>t\};\Omega)\,dt,
\label{eq:coarea}\\
\int_\Omega u\,dx
&=\int_0^\infty|\{u>t\}|\,dt.
\label{eq:layercake}
\end{align}
In particular, $\{u>t\}$ has finite relative perimeter
for almost every $t>0$.
If $V\in L^\infty(\Omega)$, Fubini's theorem also yields
$$
\int_\Omega Vu\,dx
=
\int_0^\infty\int_{\{u>t\}}V\,dx\,dt.
$$
We refer to~\cite{AFP00} for these facts.

Throughout the pairing arguments, we use the conjugate
exponents $
r:=\frac{d+1}{d}$, $q:=d+1$ and $\frac1r+\frac1q=1$.
Thus $u\,\operatorname{div}z\in L^1(\Omega)$ whenever
$u\in L^r(\Omega)$ and $\operatorname{div}z\in L^q(\Omega)$. For $u\in BV(\Omega)\cap L^r(\Omega)$ and
$z\in L^\infty(\Omega;\R^d)$ with
$\diver z\in L^q(\Omega)$, the pairing $(z,Du)$
is defined in the distributional sense by
\begin{equation}\label{eq:anzellotti}
\int_\Omega\varphi\,d(z,Du)
=
-\int_\Omega u\varphi\,\diver z\,dx
-\int_\Omega uz\cdot\nabla\varphi\,dx,
\qquad \varphi\in C_c^1(\Omega).
\end{equation} 
The distribution in \eqref{eq:anzellotti} is a finite
Radon measure satisfying $
|(z,Du)|\leq\|z\|_\infty|Du|$.
If $z$ has zero normal trace, the integration-by-parts
formula gives
\begin{equation}\label{eq:weak-pairing-bound}
\int_\Omega u\,\diver z\,dx
=
-\int_\Omega d(z,Du)
\leq\|z\|_\infty|Du|(\Omega).
\end{equation}
Here zero normal trace means that
$$
\int_\Omega \phi\,\operatorname{div}z\,dx
=
-\int_\Omega z\cdot\nabla\phi\,dx
\qquad
\text{for every }\phi\in C^1(\overline\Omega).
$$
These pairing and trace properties are recalled
from~\cite{Anz83}; see also their use in~\cite{CP19}.

For $p\geq1$, set $
c_p(x,y)=|x-y|^p$, and define the $p$-Wasserstein distance by
$$W_p^p(\mu,\nu)
=
\min_{\gamma\in\Pi(\mu,\nu)}
\int_{K\times K}|x-y|^p\,d\gamma(x,y).
$$
We write $c=c_2$, so the quadratic cost is $|x-y|^2$.
For $f\in C(K)$, define
$$
f^{c_p}(x)=\min_{y\in K}\{|x-y|^p-f(y)\}.
$$
In the quadratic case, this becomes
\begin{equation}\label{eq:ctransform}
f^c(x)=\min_{y\in K}\{|x-y|^2-f(y)\}.
\end{equation}
Iterated transforms are denoted by $f^{cc}$ and
$f^{ccc}$.
A function is called $c$-concave if it is the
$c$-transform of a function in $C(K)$; the analogous
convention applies to $c_p$. For probability measures on the compact set $K$,
optimal transport plans exist.
Moreover, narrow convergence is equivalent to
convergence in $W_p$ for each finite $p\geq1$.
Consequently, if $\mu_k\rightharpoonup\mu$ and
$\nu_k\rightharpoonup\nu$ narrowly, then $W_p(\mu_k,\nu_k)\rightarrow W_p(\mu,\nu)$.
For the quadratic cost, Kantorovich duality gives
\begin{equation}\label{eq:KD}
W_2^2(\mu,\nu)
=
\max_{f\in C(K)}
\left\{\int_K f^c\,d\mu+\int_K f\,d\nu\right\}.
\end{equation}
Thus, when $\nu=\rho\,dx$, the potential $f$ is
integrated against $\rho\,dx$, whereas $f^c$ is
integrated against $\mu$.

The $c$-transform satisfies
$f^{cc}\geq f$, $f^{ccc}=f^c$,
and
\begin{equation}\label{eq:ctransform-stability}
\|f^c-g^c\|_\infty\leq\|f-g\|_\infty,
\qquad
(f+a)^c=f^c-a
\quad\text{for }a\in\R.
\end{equation}
Since
\[
\bigl||x-y|^2-|x-y'|^2\bigr|
\leq2D|y-y'|,
\qquad x,y,y'\in K,
\]
every $c$-transform is $2D$-Lipschitz.
Replacing $f$ by $f^{cc}$ does not decrease the
objective in \eqref{eq:KD}, and subtracting
$f^{cc}(x_*)$ leaves it unchanged because both
measures have unit mass. The analogous statements hold for $c_p$, using
$$
\bigl||x-y|^p-|x-y'|^p\bigr|
\leq pD^{p-1}|y-y'|\qquad x,y,y'\in K.
$$
If $\mu$ is absolutely continuous with respect to
Lebesgue measure, the optimal quadratic transport plan
from $\mu$ to $\nu$ is unique and has the form $
\gamma=(\id,T)_\#\mu$ and $T_\#\mu=\nu$. Let $\psi$ be a $c$-concave optimal potential integrated
against $\mu$. With our cost convention,
\begin{equation}\label{eq:map-normalization}
T(x)=x-\frac12\nabla\psi(x)
\qquad\text{for }\mu\text{-a.e. }x.
\end{equation}
Indeed, the optimality relation
$\psi(x)+\psi^c(T(x))=|x-T(x)|^2$. For $c_p(x,y)=|x-y|^p$, $p>1$, absolute continuity
of the first marginal likewise guarantees a unique
optimal plan induced by a map.
We refer to~\cite[Theorem~1.17 and Remark~1.18]{San15}.

Finally, we recall the convex-concave form of
Sion's minimax theorem.

\begin{lemma}[Sion's minimax theorem {\cite{Sion58}}]
\label{L2.1}
Let $X$ and $Y$ be nonempty convex subsets of
Hausdorff topological vector spaces, and suppose
that $Y$ is compact.
Let $G:X\times Y\to\R$ satisfy the following conditions:
\begin{enumerate}
\item For every $y\in Y$, the function
$x\mapsto G(x,y)$ is convex and lower semicontinuous
on $X$.
\item For every $x\in X$, the function
$y\mapsto G(x,y)$ is concave and upper semicontinuous
on $Y$.
\end{enumerate}
Then
$$
\inf_{x\in X}\sup_{y\in Y}G(x,y)
=
\sup_{y\in Y}\inf_{x\in X}G(x,y).
$$
\end{lemma}

\section{The primal problem and its dual }\label{sec:duality}

We first establish existence for arbitrary prescribed
probability measures and uniqueness when at least one
of them is absolutely continuous.
We then derive an attained dual formulation
in terms of transport potentials and relative
perimeters.

Given an integer $n\geq 2$, let $\mu_1,\ldots,\mu_n\in\Prob(K)$, consider
\begin{equation}\label{eq:primal}
\inf_{\rho\in\A}J(\rho),
\quad
J(\rho):=
\sum_{i=1}^n W_2^2(\mu_i,\rho\,dx)+|D\rho|(\Omega).
\end{equation}

\subsection{Existence and uniqueness}
\begin{theorem}[Existence and uniqueness]\label{thm:existence}
Problem \eqref{eq:primal} admits a minimizer for arbitrary prescribed measures
in $\Prob(K)$. If at least one $\mu_i$ is absolutely continuous with
respect to Lebesgue measure, the minimizer is unique.
\end{theorem}
\begin{proof}
To see that the infimum is finite, consider the constant density $\rho_0=|\Omega|^{-1}$.
Since $\Omega$ is a nonempty bounded open set, we have
$0<|\Omega|<\infty$. Thus $\rho_0\geq 0$ and $
\int_\Omega \rho_0(y)\,dy=1$.
Moreover, $\rho_0$ is constant in $\Omega$, so its distributional
derivative in $\Omega$ vanishes. Consequently, $\rho_0\in BV(\Omega)$ and $|D\rho_0|(\Omega)=0$, which shows that $\rho_0\in\mathcal A$.
For each $i=1,\ldots,n$, the product measure
$\pi_i:=\mu_i\otimes\rho_0\,dy$ belongs to
$\Pi(\mu_i,\rho_0\,dx)$. By the definition of the quadratic
Wasserstein distance and the bound $|x-y|\leq D$ for
$x,y\in\overline{\Omega}$, we obtain
$$
\begin{aligned}
W_2^2(\mu_i,\rho_0\,dy)
&=\inf_{\pi\in\Pi(\mu_i,\rho_0\,dx)}
\int_{\overline{\Omega}\times\overline{\Omega}}
  |x-y|^2\,d\pi(x,y)\\
&\leq
  \int_{\overline{\Omega}\times\overline{\Omega}}
  |x-y|^2\,d(\mu_i\otimes\rho_0\,dy)(x,y)\leq D^2.
\end{aligned}
$$
Therefore,
$$
0\leq\inf_{\rho\in\mathcal A}J(\rho)
\leq J(\rho_0)
=\sum_{i=1}^{n}W_2^2(\mu_i,\rho_0\,dy)+|D\rho_0|(\Omega)
\leq nD^2<\infty.
$$
Let $(\rho_k)\subset\mathcal A$ be a minimizing sequence.
Since the infimum is finite, we may assume that $
\sup_k J(\rho_k)<\infty$. The Wasserstein terms are nonnegative, and hence $|D\rho_k|(\Omega)\leq J(\rho_k)$.
Moreover, the nonnegativity and unit mass of $\rho_k$ imply
$$
\sup_k\|\rho_k\|_{BV(\Omega)}
=\sup_k\bigl(
  \|\rho_k\|_{L^1(\Omega)}+|D\rho_k|(\Omega)
\bigr)<\infty.
$$
Since $\Omega$ is a bounded Lipschitz domain, the compact
embedding of $BV(\Omega)$ into $L^1(\Omega)$ yields a
subsequence, not relabeled, and a function $\rho\in BV(\Omega)$
such that $
\rho_k\rightarrow\rho$ strongly in $L^1(\Omega)$. Thus, the mass constraint
also passes to the limit. In particular, writing $\rho^-:=\max\{-\rho,0\}$ and using
$\rho_k\geq0$, we obtain
$$
\|\rho^-\|_{L^1(\Omega)}
\leq\|\rho-\rho_k\|_{L^1(\Omega)}
\rightarrow0.
$$
Thus $\rho\geq0$ almost everywhere. Furthermore, the lower semicontinuity of the total variation
under $L^1$ convergence gives
$$
|D\rho|(\Omega)
\leq\liminf_{k\to\infty}|D\rho_k|(\Omega).
$$
Therefore, $\rho\in\mathcal A$. The strong convergence $\rho_k\to\rho$ in $L^1(\Omega)$
implies that $\rho_k\,dx$ converges narrowly to $\rho\,dx$
as probability measures on $K$. By the lower semicontinuity
of the quadratic Wasserstein cost under narrow convergence,
$$
W_2^2(\mu_i,\rho\,dx)
\leq\liminf_{k\to\infty}W_2^2(\mu_i,\rho_k\,dx),
\qquad\text{for }i=1,\ldots,n.
$$
Together with the lower semicontinuity of the total variation,
this yields
$$
\begin{aligned}
J(\rho)
&=\sum_{i=1}^{n}W_2^2(\mu_i,\rho\,dx)+|D\rho|(\Omega)\\
&\leq\liminf_{k\to\infty}
\left[
\sum_{i=1}^{n}W_2^2(\mu_i,\rho_k\,dx)
+|D\rho_k|(\Omega)
\right]=\inf_{\eta\in\mathcal A}J(\eta).
\end{aligned}
$$
Therefore, $\rho$ is a minimizer.

For uniqueness, suppose $\mu_j$ is absolutely continuous. Let
$\nu_0\ne\nu_1$ be probability measures and let
$\gamma_a=(\id,T_a)_\#\mu_j$ be the optimal plans from $\mu_j$ to
$\nu_a$, $a=0,1$. For $0<s<1$, the plan
$\gamma_{s}:=(1-s)\gamma_0+s\gamma_1$ transports $\mu_j$ to
$\nu_s=(1-s)\nu_0+s\nu_1$. Hence
$$
 W_2^2(\mu_j,\nu_s)\leq
 (1-s)W_2^2(\mu_j,\nu_0)+sW_2^2(\mu_j,\nu_1).
$$
Suppose that equality holds. Then $\gamma_s$ attains
the minimum transport cost from $\mu_j$ to $\nu_s$,
and hence is an optimal plan. Since $\mu_j$ is
absolutely continuous, Brenier's theorem implies that
$\gamma_s$ is induced by a measurable map $T_s$, $\gamma_s=(\mathrm{id},T_s)_\#\mu_j$. 
By uniqueness of disintegration,
$$
(1-s)\delta_{T_0(x)}+s\delta_{T_1(x)}
=\delta_{T_s(x)}
\qquad\text{for }\mu_j\text{-almost every }x.
$$
Both coefficients $1-s$ and $s$ are strictly positive.
If $T_0(x)\neq T_1(x)$, the measure on the left assigns
positive mass to two distinct points and therefore
cannot be a Dirac mass. Consequently,
$T_0(x)=T_1(x)$ for $\mu_j$-almost every $x$. This implies that $\nu_0=(T_0)_\#\mu_j=(T_1)_\#\mu_j=\nu_1$,
contrary to our assumption. Therefore,
$$
W_2^2(\mu_j,\nu_s)
<
(1-s)W_2^2(\mu_j,\nu_0)
+sW_2^2(\mu_j,\nu_1)
\qquad\text{whenever }\nu_0\neq\nu_1.
$$
Thus the displayed inequality is strict. Since total variation is convex, $J$ is strictly convex
on the convex set $\A$, which proves uniqueness.
\end{proof}

\begin{remark}
The uniqueness argument is the classical transport argument underlying
Wasserstein barycenter uniqueness~\cite{bib1}. It is independent of the
dual formulation. Without an absolutely continuous prescribed measure, the theorem asserts
existence but makes no general uniqueness claim.
\end{remark}

\subsection{A dual formulation through relative perimeter}
We now seek a dual expression for the minimum in
\eqref{eq:primal}.
The mass and nonnegativity constraints lead naturally
to an inequality involving the sum of the transport
potentials.
Before establishing the dual formula, we explain
why imposing pointwise cancellation can give
a smaller value than the primal minimum. To motivate a dual formulation, let
$f_1,\ldots,f_n\in C(K)$ and
$z\in C_c^1(\Omega;\mathbb R^d)$ with
$\|z\|_\infty\leq1$.
For every $\rho\in\mathcal A$, Kantorovich duality gives
$$
W_2^2(\mu_i,\rho\,dx)
\geq
\int_K f_i^c\,d\mu_i+\int_\Omega f_i\rho\,dx,
\qquad i=1,\ldots,n.
$$
Moreover, the definition of the total variation yields $
|D\rho|(\Omega)
\geq\int_\Omega\rho\,\operatorname{div}z\,dx$. Adding these inequalities, we obtain
$$
J(\rho)
\geq
\sum_{i=1}^{n}\int_K f_i^c\,d\mu_i
+\int_\Omega
\left(\sum_{i=1}^{n}f_i+\operatorname{div}z\right)\rho\,dx.
$$
Since every $\rho\in\mathcal A$ is nonnegative and has
unit mass, the condition $
\sum_{i=1}^{n}f_i+\operatorname{div}z\geq\lambda$ in $\Omega$
implies
$$
J(\rho)\geq
\lambda+\sum_{i=1}^{n}\int_K f_i^c\,d\mu_i.
$$
Thus a lower bound independent of $\rho$ can be obtained
without requiring
$\sum_{i=1}^{n}f_i+\operatorname{div}z$ to vanish identically.
Imposing this identity gives a more restrictive class
of competitors. The following example shows that the
resulting supremum can be strictly smaller than the
primal minimum.

\begin{example}\label{ex:boundary}
Let $\Omega=(0,1)$, $n=2$, and $
\mu_1=\mu_2=10\mathbf{1}_{(0,1/10)}\,dx$. In this case, the candidate dual value is
$$
\mathsf D_0
=
\sup\left\{
\sum_{i=1}^{2}\int_{[0,1]}f_i^c\,d\mu_i:
\begin{array}{l}
f_i\in C([0,1]),\
z\in C_c^1((0,1)),\\
\|z\|_\infty\leq1,\
f_1+f_2=-z'
\text{ in }(0,1)
\end{array}
\right\}.
$$
Since $z$ is compactly supported in $(0,1)$, its derivative
vanishes in a neighborhood of $0$. The constraint and
continuity of the potentials therefore give $f_1(0)+f_2(0)=0$.
Taking $y=0$ in the definition
$$
f_i^c(x)
=\min_{y\in[0,1]}\bigl\{|x-y|^2-f_i(y)\bigr\},
$$
we obtain that $f_i^c(x)\leq x^2-f_i(0)$. Since $\mu_1=\mu_2$, every admissible choice satisfies
$$
\begin{aligned}
\sum_{i=1}^{2}\int_{[0,1]}f_i^c\,d\mu_i
&\leq
10\int_0^{1/10}
\bigl(2x^2-f_1(0)-f_2(0)\bigr)\,dx\\
&=20\int_0^{1/10}x^2\,dx
=\frac1{150}.
\end{aligned}
$$
Hence $\mathsf D_0\leq1/150$. We next derive a strictly larger lower bound for the
primal problem. Let $\rho\in\mathcal A$ and set $
M:=|D\rho|((0,1))$. If $M\leq1$, the one-dimensional $BV$ oscillation estimate,
together with $\int_0^1\rho\,dx=1$, yields
$$
\operatorname*{ess\,sup}_{(0,1)}\rho
\leq
\operatorname*{ess\,inf}_{(0,1)}\rho+M
\leq1+M\leq2.
$$
It follows that
$$
\int_{1/5}^{1}\rho(x)\,dx
=1-\int_0^{1/5}\rho(x)\,dx
\geq1-\frac25=\frac35.
$$
For each prescribed measure, all the source mass lies in
$[0,1/10]$. Any transport plan to $\rho\,dx$ must therefore
move at least $3/5$ of its mass over a distance of at
least $1/10$. Consequently,
$$
W_2^2(\mu_i,\rho\,dx)
\geq\frac35\left(\frac1{10}\right)^2,
\qquad i=1,2,
$$
and thus
$$
J(\rho)\geq
2\cdot\frac35\left(\frac1{10}\right)^2
=\frac3{250}.
$$
If $M>1$, then $J(\rho)\geq M>1$, so the same lower bound
holds. We conclude that
$$
\mathsf D_0
\leq\frac1{150}
<\frac3{250}
\leq\inf_{\rho\in\mathcal A}J(\rho).
$$
Thus the supremum in the proposed dual formulation is
strictly smaller than the primal minimum. Restricting the potentials further to $C_c(\Omega)$
can only decrease this supremum, so the strict duality
gap persists under this additional restriction.
\end{example}

This example concerns the simultaneous imposition of compact support on the vector field
and the pointwise cancellation condition $\sum_i f_i+\operatorname{div}z=0$.
It does not apply to the weak vector-field formulation
of Carlier and Poon~\cite[Theorem~3.3]{CP19},
which uses zero normal trace and a nonnegative multiplier.
To retain the mass and nonnegativity constraints explicitly, for $V\in C(\overline\Omega)$, we introduce the auxiliary minimization problem
\begin{equation}\label{eq:mV}
m(V)=
\inf_{\rho\in\A}
\left\{
|D\rho|(\Omega)+\int_\Omega V\rho\,dx
\right\}.
\end{equation}
The following lemma expresses this value in terms
of relative perimeters and provides the basis
for the dual formulation.

\begin{lemma}\label{lem:coarea-dual}
The infimum in \eqref{eq:mV} is attained, and
\begin{equation}\label{eq:set-ratio}
 m(V)=\inf_{\substack{E\subset\Omega,\ |E|>0\\
                         \Per(E;\Omega)<\infty}}
 \frac{\Per(E;\Omega)+\int_E V\dd x}{|E|}.
\end{equation}
Moreover, $m$ is concave and
$|m(V)-m(W)|\leq\|V-W\|_\infty$.
\end{lemma}
\begin{proof}
Since $V\in C(K)$ and $K$ is compact, $V$ is bounded.
For every $\rho\in\mathcal A$, nonnegativity and unit
mass imply $\left|\int_\Omega V\rho\,dx\right|
\leq\|V\|_\infty$.
In particular, the functional is bounded below by
$-\|V\|_\infty$. Testing with the constant density
$\rho_0=|\Omega|^{-1}$ gives
$$
-\|V\|_\infty
\leq m(V)
\leq\frac1{|\Omega|}\int_\Omega V\,dx
\leq\|V\|_\infty.
$$
Let $(\rho_k)\subset\mathcal A$ be a minimizing sequence.
After discarding finitely many terms, we may assume that
$|D\rho_k|(\Omega)+\int_\Omega V\rho_k\,dx
\leq m(V)+1$. Consequently,
$|D\rho_k|(\Omega)\leq2\|V\|_\infty+1$.
Together with $\|\rho_k\|_{L^1(\Omega)}=1$, this shows
that $(\rho_k)$ is bounded in $BV(\Omega)$. The compactness theorem for $BV$ functions therefore
yields a subsequence converging strongly in $L^1(\Omega)$
to some $\rho\in\mathcal A$, as in the preceding
existence proof. Since
$$
\left|\int_\Omega V(\rho_k-\rho)\,dx\right|
\leq\|V\|_\infty\|\rho_k-\rho\|_{L^1(\Omega)}
\rightarrow0,
$$
the lower semicontinuity of the total variation gives
$$
|D\rho|(\Omega)+\int_\Omega V\rho\,dx
\leq
\liminf_{k\to\infty}
\left(
|D\rho_k|(\Omega)+\int_\Omega V\rho_k\,dx
\right)
=m(V).
$$
Since $\rho\in\mathcal A$, equality holds, and the
infimum defining $m(V)$ is attained. 

Set
$$
a:=
\inf_{\substack{E\subset\Omega,\ |E|>0\\
\operatorname{Per}(E;\Omega)<\infty}}
\frac{\operatorname{Per}(E;\Omega)+\int_E V\,dx}{|E|}.
$$
Since the relative perimeter is nonnegative and $\operatorname{Per}(\Omega;\Omega)=0$. Taking
$E=\Omega$ gives $\min_K V\leq a
\leq\frac1{|\Omega|}\int_\Omega V\,dx$. In particular, $a$ is finite. For every set $E$ in the above infimum, the normalized
characteristic function $
\rho_E:=\frac{\mathbf1_E}{|E|}$ belongs to $\mathcal A$. By homogeneity of the total
variation,
$$
|D\rho_E|(\Omega)+\int_\Omega V\rho_E\,dx
=
\frac{\operatorname{Per}(E;\Omega)+\int_E V\,dx}{|E|}.$$
Testing the definition of $m(V)$ with $\rho_E$ and then
taking the infimum over $E$ yields $m(V)\leq a$. To prove the reverse inequality, fix $\rho\in\mathcal A$
and let $E_t:=\{x\in\Omega:\rho(x)>t\}$, $t>0$.
Since $\rho\in BV(\Omega)$,  coarea formula gives $\operatorname{Per}(E_t;\Omega)<\infty$ for almost every $t>0$.
The layer-cake representation and the mass constraint give
$\int_0^\infty|E_t|\,dt=1$.
Since $\int_\Omega |V|\rho\,dx\leq\|V\|_\infty<\infty$,
Fubini's theorem and the layer-cake representation yield
$$
\int_\Omega V\rho\,dx
=\int_0^\infty\int_{E_t}V\,dx\,dt.
$$
For almost every $t$ with $|E_t|>0$, the definition of $a$
gives $\operatorname{Per}(E_t;\Omega)+\int_{E_t}V\,dx
\geq a|E_t|$. The same inequality holds when $|E_t|=0$, since then
$\mathbf1_{E_t}=0$ almost everywhere and its relative
perimeter is zero. Integrating in $t$, we obtain
$$
\begin{aligned}
|D\rho|(\Omega)+\int_\Omega V\rho\,dx
&=\int_0^\infty
\left(
\operatorname{Per}(E_t;\Omega)+\int_{E_t}V\,dx
\right)\,dt\\
&\geq a\int_0^\infty|E_t|\,dt
=a.
\end{aligned}
$$
Taking the infimum over $\rho\in\mathcal A$ proves
$m(V)\geq a$, and hence the claimed formula. Concavity
follows because \eqref{eq:mV} is an infimum of affine functions of $V$. Finally, nonnegativity and unit mass imply
$$
\left|\int_\Omega(V-W)\rho\,dx\right|
\leq\|V-W\|_\infty
\qquad\text{for every }\rho\in\mathcal A.
$$
It follows by taking infima that $
m(V)\leq m(W)+\|V-W\|_\infty$.
Interchanging $V$ and $W$ yields $
|m(V)-m(W)|\leq\|V-W\|_\infty$.
\end{proof}

Fix $x_*\in K$ and introduce the compact convex set
\begin{equation}\label{eq:normalized-potentials}
 \K=\{f\in C(K):f(x_*)=0,\ \Lip(f)\leq2D\}.
\end{equation}
Compactness is with respect to the uniform norm. Notice that $\|f\|_\infty\leq2D^2$ for $f\in\K$.
For $\boldsymbol f=(f_1,\ldots,f_n)\in C(K)^n$, write
$$
V_{\boldsymbol f}:=\sum_{i=1}^{n}f_i,
\qquad
B(\boldsymbol f):=\sum_{i=1}^{n}\int_K f_i^c\,d\mu_i.
$$

\begin{theorem}\label{thm:duality}
For arbitrary $\mu_1,\ldots,\mu_n\in\Prob(K)$,
\begin{equation}\label{eq:dual-m}
 \min_{\rho\in\A}J(\rho)
 =\max_{\boldsymbol f\in\K^n}
       \left\{B(\boldsymbol f)+m(V_{\boldsymbol f})\right\}.
\end{equation}
Equivalently, the minimum equals
\begin{equation}\label{eq:dual-sets}
 \max_{\substack{f_1,\ldots,f_n\in C(K),\ \lambda\in\R\\
 \lambda|E|\leq\Per(E;\Omega)+\int_E\sum_i f_i\dd x\
 \text{ for all finite-perimeter }E\subset\Omega}}
 \left\{\lambda+\sum_{i=1}^n\int_K f_i^c\dd\mu_i\right\}.
\end{equation}
The maximum in \eqref{eq:dual-sets} can be attained with
$f_i\in\K$ and $\lambda=m(V_{\boldsymbol f})$.
\end{theorem}
\begin{proof}
For $\rho\in\mathcal A$ and
$\boldsymbol f=(f_1,\ldots,f_n)\in\mathcal K^n$, define
$$
L(\rho,\boldsymbol f)
:=
|D\rho|(\Omega)
+\int_\Omega V_{\boldsymbol f}\rho\,dx
+B(\boldsymbol f).$$
We first justify restricting the Kantorovich dual problem
to $\mathcal K$. For $\rho\in\mathcal A$, we regard $\nu:=\rho\,dx$ as
a probability measure on $K$ by setting $
\nu(A):=\int_{A\cap\Omega}\rho(x)\,dx$
for every Borel set $A\subset K$. In particular, for $f\in C(K)$, $\nu(\partial\Omega)=0$ and $
\int_K f\,d\nu=\int_\Omega f\rho\,dx$ . For each $i$, Kantorovich duality gives
$$
W_2^2(\mu_i,\nu)
=
\sup_{f\in C(K)}
\left\{
\int_K f^c\,d\mu_i+\int_K f\,d\nu
\right\}.$$
Given $f\in C(K)$, its double $c$-transform satisfies $
f^{cc}\geq f$ and $(f^{cc})^c=f^c$. Consequently, replacing $f$ by $f^{cc}$ does not decrease
the dual objective$$
\int_K (f^{cc})^c\,d\mu_i+\int_K f^{cc}\,d\nu
\geq
\int_K f^c\,d\mu_i+\int_K f\,d\nu.$$
Since
$$
f^{cc}(y)
=\inf_{x\in K}\bigl\{|x-y|^2-f^c(x)\bigr\}$$
and
$$
\bigl||x-y|^2-|x-y'|^2\bigr|
\leq2D|y-y'|,
\qquad x,y,y'\in K,$$
the function $f^{cc}$ is $2D$-Lipschitz. The Kantorovich objective is invariant under the normalization
$\widetilde f:=f-f(x_*)$, since $
\widetilde f^{\,c}=f^c+f(x_*)$
and both marginal measures have unit mass. Thus optimal potentials can be chosen
in $\mathcal K$, and Kantorovich duality yields
$$
J(\rho)
=\max_{\boldsymbol f\in\mathcal K^n}
L(\rho,\boldsymbol f)
\qquad\text{for every }\rho\in\mathcal A.$$
We next verify the hypotheses of Lemma \ref{L2.1}.
Equip $\mathcal A$ with the topology inherited from
$L^1(\Omega)$ and $\mathcal K^n$ with the product uniform
topology. Both sets are convex, and $\mathcal K^n$ is
compact. For fixed $\boldsymbol f$, the map
$\rho\mapsto L(\rho,\boldsymbol f)$ is convex and lower
semicontinuous,
while the linear term is continuous because
$V_{\boldsymbol f}$ is bounded. For fixed $\rho$, the map
$\boldsymbol f\mapsto L(\rho,\boldsymbol f)$ is concave.
Indeed, for $f,g\in C(K)$ and $s\in[0,1]$,
$$
\begin{aligned}
\bigl((1-s)f+sg\bigr)^c(x)
&=\min_{y\in K}
\left\{
(1-s)\bigl(|x-y|^2-f(y)\bigr)
+s\bigl(|x-y|^2-g(y)\bigr)
\right\}\\
&\geq(1-s)f^c(x)+sg^c(x).
\end{aligned}
$$
It is also continuous, since $
\|f^c-g^c\|_\infty\leq\|f-g\|_\infty$ and $\int_\Omega\rho\,dx=1$. Therefore, we have
$$
\inf_{\rho\in\mathcal A}J(\rho)
=\inf_{\rho\in\mathcal A}
  \sup_{\boldsymbol f\in\mathcal K^n}
  L(\rho,\boldsymbol f)=\sup_{\boldsymbol f\in\mathcal K^n}
  \inf_{\rho\in\mathcal A}
  L(\rho,\boldsymbol f)=\sup_{\boldsymbol f\in\mathcal K^n}
  \bigl\{B(\boldsymbol f)+m(V_{\boldsymbol f})\bigr\}.
$$
The last supremum is attained. Indeed, the stability
of the $c$-transform and the Lipschitz continuity of $m$
established in Lemma~\ref{lem:coarea-dual} imply
$$
\begin{aligned}
|B(\boldsymbol f)-B(\boldsymbol g)|
&\leq\sum_{i=1}^{n}\|f_i-g_i\|_\infty,\\
|m(V_{\boldsymbol f})-m(V_{\boldsymbol g})|
&\leq\|V_{\boldsymbol f}-V_{\boldsymbol g}\|_\infty
\leq\sum_{i=1}^{n}\|f_i-g_i\|_\infty.
\end{aligned}
$$
Thus the objective is continuous on the compact set
$\mathcal K^n$. Together with primal attainment, this proves
\eqref{eq:dual-m}. It remains to prove \eqref{eq:dual-sets}. By Lemma~\ref{lem:coarea-dual}, the constraints $\lambda|E|
\leq\operatorname{Per}(E;\Omega)
+\int_E V_{\boldsymbol f}\,dx$
for all measurable sets $E\subset\Omega$ of positive
measure and finite relative perimeter are equivalent to $
\lambda\leq m(V_{\boldsymbol f})$. Hence, for any admissible pair
$(\boldsymbol f,\lambda)\in C(K)^n\times\mathbb R$
and any $\rho\in\mathcal A$,
$$
\lambda+B(\boldsymbol f)\leq m(V_{\boldsymbol f})+B(\boldsymbol f)\leq |D\rho|(\Omega)
     +\int_\Omega V_{\boldsymbol f}\rho\,dx
     +B(\boldsymbol f)\leq J(\rho),
$$
where the last inequality follows from Kantorovich
duality. Finally, let $\boldsymbol f^*\in\mathcal K^n$ attain
the maximum in \eqref{eq:dual-m} and set
$\lambda^*:=m(V_{\boldsymbol f^*})$.
The pair $(\boldsymbol f^*,\lambda^*)$ satisfies all
the perimeter constraints and $
\lambda^*+B(\boldsymbol f^*)
=\min_{\rho\in\mathcal A}J(\rho)$.
This proves the stated dual equality \eqref{eq:dual-sets}.
\end{proof}

\begin{remark}
The dual constraints and objective are invariant under
\[
f_i\mapsto f_i+a_i,
\qquad
\lambda\mapsto\lambda+\sum_{i=1}^{n}a_i,
\qquad a_i\in\mathbb R.
\]
Thus the normalization $f_i(x_*)=0$ only fixes the
additive constants and imposes no boundary condition.
Moreover, the relative perimeter
$\operatorname{Per}(E;\Omega)$ counts only the perimeter
inside $\Omega$, consistently with the total variation
$|D\rho|(\Omega)$ in the primal problem.
\end{remark}

\section{Optimality conditions and transport maps}\label{sec:optimality}
The attained dual formulation provides a direct
characterization of the minimizers.
For a fixed dual maximizer, optimality is equivalent
to equality in each transport inequality and in
the total variation bound.
We first state these conditions and their consequences
for superlevel sets, and then derive a transport-map
representation and a criterion for the uniform
density to be optimal.

\begin{theorem}[Transport and TV optimality]\label{thm:optimality}
Let $(\boldsymbol f,\lambda)$ be a maximizer of
\eqref{eq:dual-sets}. A density $\rho\in\A$ minimizes $J$ if and only if
\begin{align}
 &W_2^2(\mu_i,\rho\,dx)
 =\int_Kf_i^c\,d\mu_i+\int_\Omega f_i\rho\,dx,\qquad i=1,\ldots,n,\label{eq:OT-equality}\\
 &|D\rho|(\Omega)+\int_\Omega V_{\boldsymbol f}\rho\,d x=\lambda.\label{eq:TV-equality}
\end{align}
For such a density, almost every $t>0$ with $|E_t|>0$, where
$E_t=\{\rho>t\}$, satisfies
\begin{equation}\label{eq:levelsets}
 \Per(E_t;\Omega)+\int_{E_t}V_{\boldsymbol f}\dd x
      =\lambda|E_t|.
\end{equation}
Thus these superlevel sets minimize the ratio in
\eqref{eq:set-ratio} for the same potential $V_{\boldsymbol f}$.
\end{theorem}
\begin{proof}
Fix a dual maximizer $(\boldsymbol f,\lambda)$. By the strong duality, $
\lambda+B(\boldsymbol f)=\min_{\rho\in\mathcal A}J(\rho)$. For any $\rho\in\mathcal A$, Kantorovich duality gives
$$
W_2^2(\mu_i,\rho\,dx)
-\int_K f_i^c\,d\mu_i
-\int_\Omega f_i\rho\,dx
\geq0,
\qquad i=1,\ldots,n.
$$
Since $(\boldsymbol f,\lambda)$ is a dual maximizer,
we have $\lambda=m(V_{\boldsymbol{f}})$. Consequently, by the definition of $m(V)$,
$$
|D\rho|(\Omega)+\int_\Omega V_{\boldsymbol{f}}\rho\,dx-\lambda
=
|D\rho|(\Omega)+\int_\Omega V_{\boldsymbol{f}}\rho\,dx-m(V_{\boldsymbol{f}})
\geq0
$$
for every $\rho\in\mathcal A$. Thus the identity
$$
\begin{aligned}
J(\rho)-(\lambda+B({\boldsymbol{f}}))
&=\sum_{i=1}^{n}
\left[
W_2^2(\mu_i,\rho\,dx)
-\int_K f_i^c\,d\mu_i
-\int_\Omega f_i\rho\,dx
\right]\\
&\quad+
\left[
|D\rho|(\Omega)+\int_\Omega V_{\boldsymbol{f}}\rho\,dx-\lambda
\right]
\end{aligned}
$$
expresses the primal--dual gap as a sum of nonnegative
terms. Since $\lambda+B({\boldsymbol{f}})$ is the primal minimum,
$\rho$ is a minimizer if and only if this sum vanishes.
This occurs precisely when every term is zero, giving
\eqref{eq:OT-equality} and \eqref{eq:TV-equality}. Now let $\rho$ be a minimizer and set
$E_t:=\{\rho>t\}$ for $t>0$.
By the coarea formula, $E_t$ has finite relative
perimeter for almost every $t>0$. Dual feasibility gives
$$
G(t):=\operatorname{Per}(E_t;\Omega)
+\int_{E_t}V_{\boldsymbol{f}}\,dx-\lambda|E_t|\geq0
$$
for almost every such $t$; when $|E_t|=0$, all three
terms vanish.
Using the coarea and layer-cake formulas, we obtain
$$
\begin{aligned}
\int_0^\infty G(t)\,dt
&=|D\rho|(\Omega)+\int_\Omega V_{\boldsymbol{f}}\rho\,dx
  -\lambda\int_\Omega\rho\,dx\\
&=|D\rho|(\Omega)+\int_\Omega V_{\boldsymbol{f}}\rho\,dx-\lambda
=0,
\end{aligned}
$$
where the last equality is \eqref{eq:TV-equality}.
The term involving $V$ is integrable because
$\int_\Omega|V_{\boldsymbol{f}}|\rho\,dx\leq\|V_{\boldsymbol{f}}\|_\infty$.
Since $G\geq0$ almost everywhere and its integral is zero,
$G(t)=0$ for almost every $t>0$, which proves
\eqref{eq:levelsets}. Since $\lambda=m(V_{\boldsymbol{f}})$, the level-set identity
\eqref{eq:levelsets} gives
$$
\frac{\operatorname{Per}(E_t;\Omega)+\int_{E_t}V_{\boldsymbol{f}}\,dx}
     {|E_t|}
=m(V_{\boldsymbol{f}})
$$
for almost every $t>0$ with $|E_t|>0$.
By Lemma~\ref{lem:coarea-dual}, $m(V)$ is the infimum
of this ratio over all measurable subsets of $\Omega$
with positive measure and finite relative perimeter.
Hence each such superlevel set $E_t$ attains this infimum.
\end{proof}

When at least one of the prescribed measures is
absolutely continuous, the transport optimality
conditions yield the following representation
of the unique minimizing measure.
\begin{corollary}[Map representation]\label{cor:map}
Suppose that $\mu_j$ is absolutely continuous
for some $j\in\{1,\ldots,n\}$, and fix any
dual maximizer $(\boldsymbol f,\lambda)$. The unique
minimizer satisfies
\begin{equation}\label{eq:pushforward}
 \rho\,dx=(T_j)_\#\mu_j,
 \qquad T_j(x)=x-\frac12\nabla f_j^c(x)\text{ for }\mu_j-a.e. x.
\end{equation}
Equivalently, setting
$\Phi_j(x)=\tfrac12(|x|^2-f_j^c(x))$, one has
$\rho\dd x=(\nabla\Phi_j)_\#\mu_j$, and $\Phi_j$ is convex.
\end{corollary}
\begin{proof}
Equation~\eqref{eq:OT-equality} makes $(f_j^c,f_j)$ an optimal
Kantorovich pair. The transport theorem and
\eqref{eq:map-normalization} yield \eqref{eq:pushforward}.
For the convexity assertion, extend the transform to $\R^d$ by the
same infimum over $y\in K$ and observe that
$$
 \Phi_j(x)=\sup_{y\in K}\left\{x\cdot y-\frac{|y|^2}{2}+\frac{f_j(y)}2\right\}.
$$
This is a supremum of affine functions.
\end{proof}

For the uniform density, the total variation vanishes,
and the optimality conditions reduce to the following
criterion.
\begin{corollary}[A criterion for the uniform density]\label{cor:uniform}
The density $\rho_0=|\Omega|^{-1}$ minimizes $J$ if and only if there
exist potentials $f_1,\ldots,f_n\in C(K)$ such that $(f_i^c,f_i)$ is
optimal for $(\mu_i,\rho_0\dd x)$ for each $i$, and, writing
$V=\sum_i f_i$ and $\overline V=|\Omega|^{-1}\int_\Omega V\dd x$,
\begin{equation}\label{eq:uniform-criterion}
 \Per(E;\Omega)+\int_E(V-\overline V)\dd x\geq0
 \quad\text{for every finite-perimeter }E\subset\Omega.
\end{equation}
\end{corollary}
\begin{proof}
If $\rho_0$ minimizes $J$, choose a dual maximizer and apply
Theorem~\ref{thm:optimality}. Since $\TV(\rho_0)=0$, it gives
$\lambda=\overline V$, and the dual constraints give
\eqref{eq:uniform-criterion}. Conversely, the stated conditions make
$(\boldsymbol f,\overline V)$ dual feasible and make its value equal
to $J(\rho_0)$. Weak duality then proves optimality.
\end{proof}

\section{Vector-field duality via entropic approximation}\label{sec5}
Following Carlier and Poon~\cite[Section~3]{CP19},
we use an entropic approximation to establish an
attained vector-field dual.
The positivity of the regularized minimizers allows
us to differentiate the transport costs with respect
to the density, even when the prescribed measures
are singular.
The resulting optimality conditions involve the sum
of the Kantorovich potentials and an explicit
multiplier for the mass constraint.
We pass to the limit using uniform estimates
and then relate the resulting formulation
to the perimeter dual.

\subsection{Vector fields and total variation}

We introduce the vector fields used to represent
the total variation and recall the associated
integration-by-parts formula.
The lemma below will allow us to recover such
a field from a variational inequality. The embedding $BV(\Omega)\hookrightarrow L^r(\Omega)$ is compact,
including when $d=1$. For $z\in L^\infty(\Omega;\R^d)$ with
$\diver z=\xi\in L^q(\Omega)$, we say that its normal trace is zero if
\begin{equation}\label{eq:weak-zero-trace}
 \int_\Omega \xi u\dd x=-\int_\Omega z\cdot\nabla u\dd x
 \quad\text{for every }u\in C^1(\overline\Omega).
\end{equation}
This identity includes both the distributional divergence relation
and the zero boundary trace. Write
$$
 \Gamma_q=\left\{\xi\in L^q(\Omega):
 \begin{array}{l}
 \xi=\diver z\text{ for some }z\in L^\infty(\Omega;\R^d),\\
 \|z\|_\infty\leq1\text{ and \eqref{eq:weak-zero-trace} holds}
 \end{array}\right\}.
$$

\begin{lemma}\label{lem:TV-polar}
The set $\Gamma_q$ is convex and norm-closed in $L^q(\Omega)$. With
total variation extended by $+\infty$ outside $BV(\Omega)$,
\begin{equation}\label{eq:polar-support}
 |Du|(\Omega)=\sup_{\xi\in\Gamma_q}\int_\Omega\xi u\dd x
 \quad\text{for }u\in L^r(\Omega).
\end{equation}
If $\xi\in L^q(\Omega)$ satisfies $
\int_\Omega \xi u\,dx\leq|Du|(\Omega)$ for every $u\in BV(\Omega)\cap L^\infty(\Omega)$, then $\xi\in\Gamma_q$.
\end{lemma}
\begin{proof}
Convexity is immediate. If $\xi_k\to\xi$ in $L^q$ and $z_k$ are
corresponding fields, weak-* compactness gives $z_k\rightharpoonup^*z$
in $L^\infty$ along a subsequence, with $\|z\|_\infty\leq1$.
Passing to the limit in \eqref{eq:weak-zero-trace} proves closedness.
The upper bound in \eqref{eq:polar-support} follows from
\eqref{eq:weak-pairing-bound}; the opposite bound follows from
\eqref{eq:TV}, since compactly supported smooth fields are allowed.
If $u\in L^r(\Omega)\setminus BV(\Omega)$, then
$u\in L^1(\Omega)$ because $\Omega$ is bounded.
By the definition of $BV(\Omega)$,
$$
\sup_{\substack{z\in C_c^1(\Omega;\mathbb R^d)\\
                 \|z\|_\infty\leq1}}
\int_\Omega u\,\operatorname{div}z\,dx
=+\infty.
$$
Since $\operatorname{div}z\in\Gamma_q$ for every such
vector field, the supremum in \eqref{eq:polar-support}
is also infinite. This proves \eqref{eq:polar-support}
for all $u\in L^r(\Omega)$. Finally, suppose that $\xi\in L^q(\Omega)$ satisfies $
\int_\Omega\xi u\,dx\leq|Du|(\Omega)$ for every $u\in BV(\Omega)\cap L^\infty(\Omega)$. Let $u\in BV(\Omega)\cap L^r(\Omega)$, where $1\leq r<\infty$.
For $M>0$, define $
T_M(s):=\max\{-M,\min\{s,M\}\}$ and $u_M:=T_M(u)$.
Then $|u_M|\leq M$ almost everywhere, so
$u_M\in L^\infty(\Omega)$. We next verify that
$u_M\in BV(\Omega)$ and that truncation does not increase
the total variation. By smooth approximation in $BV(\Omega)$, there exists
a sequence $v_k\in C^\infty(\Omega)\cap W^{1,1}(\Omega)$
such that $
v_k\rightarrow u$ in $L^1(\Omega)$ and $
\int_\Omega|\nabla v_k|\,dx
\rightarrow |Du|(\Omega)$.
Since $T_M$ is $1$-Lipschitz,
$$
\|T_M(v_k)-u_M\|_{L^1(\Omega)}
\leq\|v_k-u\|_{L^1(\Omega)}
\rightarrow0.$$
The Sobolev chain rule also gives
$T_M(v_k)\in W^{1,1}(\Omega)$ and $
|\nabla T_M(v_k)|
\leq|\nabla v_k|$ almost everywhere in $\Omega$. Therefore, by the lower semicontinuity of the total
variation under $L^1$ convergence,
$$
\begin{aligned}
|Du_M|(\Omega)
&\leq\liminf_{k\to\infty}
      \int_\Omega|\nabla T_M(v_k)|\,dx\\
&\leq\lim_{k\to\infty}
      \int_\Omega|\nabla v_k|\,dx
=|Du|(\Omega)<\infty.
\end{aligned}
$$
Thus $u_M\in BV(\Omega)\cap L^\infty(\Omega)$ and
$|Du_M|(\Omega)\leq|Du|(\Omega)$. Finally, the definition of $T_M$ yields
$$
|u-u_M|=
\begin{cases}
0, & |u|\leq M,\\
|u|-M, & |u|>M.
\end{cases}
$$
Consequently,
$$
\|u-u_M\|_{L^r(\Omega)}^r
=
\int_{\{|u|>M\}}(|u|-M)^r\,dx
\leq
\int_\Omega |u|^r\mathbf1_{\{|u|>M\}}\,dx.
$$
As $M\to\infty$, the integrand on the right converges
to zero almost everywhere and is bounded by
$|u|^r\in L^1(\Omega)$. The dominated convergence theorem
therefore implies $
u_M\rightarrow u$ strongly in $L^r(\Omega)$. 
Thus $
\int_\Omega\xi u\,dx\leq|Du|(\Omega)$.
For $u\in L^r(\Omega)\setminus BV(\Omega)$, this inequality
holds trivially because $|Du|(\Omega)=+\infty$.
It therefore holds for every $u\in L^r(\Omega)$. If $\xi\notin\Gamma_q$, the Hahn-Banach separation theorem,
applied to the nonempty closed convex set
$\Gamma_q\subset L^q(\Omega)$, yields some
$u\in L^r(\Omega)=(L^q(\Omega))^*$ such that
$$
\int_\Omega\xi u\,dx
>
\sup_{\eta\in\Gamma_q}\int_\Omega\eta u\,dx
=|Du|(\Omega),
$$
contradicting the preceding inequality.
Hence $\xi\in\Gamma_q$.
\end{proof}

\subsection{Entropy regularization and estimates}

We approximate $J$ by adding an entropy term.
The resulting minimizers are bounded away from zero
for each fixed regularization parameter.
The estimates established below also provide the
compactness needed to pass to the limit. Write
$$
F(\rho)=\sum_{i=1}^nW_2^2(\mu_i,\rho\,d x),\qquad
 H(\rho)=\int_\Omega\rho\log\rho\,d x,
 \qquad J_\varepsilon(\rho)
:=F(\rho)+|D\rho|(\Omega)+\varepsilon H(\rho),
$$
where $0\log0=0$ and $\varepsilon>0$.
Since $BV(\Omega)\subset L^r(\Omega)$ with $r>1$, $H$ is finite on
$\A$. We shall use the elementary transport estimate
\begin{equation}\label{eq:transport-L1}
 |F(\eta)-F(\rho)|\leq 2nD^2\|\eta-\rho\|_{L^1(\Omega)},
 \qquad \rho,\eta\in\A.
\end{equation}
Indeed, for each $i$, choose an optimal potential
$f_i\in\mathcal K$ for $(\mu_i,\eta\,dx)$.
Using the same potential as a dual competitor for
$(\mu_i,\rho\,dx)$ and subtracting gives
$$
W_2^2(\mu_i,\eta\,dx)-W_2^2(\mu_i,\rho\,dx)
\leq\int_\Omega f_i(\eta-\rho)\,dx
\leq2D^2\|\eta-\rho\|_{L^1(\Omega)},
$$
where we used
$|f_i(x)|=|f_i(x)-f_i(x_*)|\leq2D|x-x_*|\leq2D^2$.
Summing over $i$ yields the corresponding upper bound
for $F(\eta)-F(\rho)$. Interchanging $\rho$ and $\eta$ gives the reverse
bound and proves the stated estimate.

\begin{lemma}\label{lem:entropy-bounds}
For every $\varepsilon>0$, $J_\varepsilon$ has a unique minimizer
$\rho_\varepsilon\in\A$, $
|D\rho_\varepsilon|(\Omega)\leq nD^2$, and $\infess_\Omega\rho_\varepsilon>0$. Along a sequence $\varepsilon\to 0$ there is
$\rho_*\in\A$ such that $\rho_\varepsilon\to\rho_*$ strongly in
$L^r$, $\rho_*$ minimizes $J$, and
\begin{equation}\label{eq:entropy-log-bounds}
 \begin{aligned}
 b_\varepsilon^-&:=-\varepsilon\log\min\{1,\rho_\varepsilon\}
       &&\text{is bounded in }L^\infty,\\
 b_\varepsilon^+&:=\varepsilon\log\max\{1,\rho_\varepsilon\}
       &&\rightarrow0\text{ in every finite }L^s.
 \end{aligned}
\end{equation}
Here $\infess$ denotes the essential infimum. No positive lower bound
for $\rho_\varepsilon$ independent of $\varepsilon$ is asserted.
\end{lemma}
\begin{proof}
Our proof is divided into three steps.

\noindent
\textit{Step 1.} By Jensen's inequality, $
H(\rho)\geq-\log|\Omega|$ for every $\rho\in\mathcal A$. Hence $J_\varepsilon(\rho)
\geq |D\rho|(\Omega)-\varepsilon\log|\Omega|$. A minimizing sequence is therefore bounded in $BV(\Omega)$.
The compact embedding into $L^1(\Omega)$, the lower
semicontinuity of the total variation and entropy, and
the continuity of $F$ under $L^1$ convergence yield a
minimizer $\rho_\varepsilon\in\mathcal A$.
Uniqueness follows from the strict convexity of $H$
and the convexity of the remaining terms. Let $\rho_0=|\Omega|^{-1}$. Since
$H(\rho_0)=-\log|\Omega|$, comparison with $\rho_0$ gives
$$
\begin{aligned}
|D\rho_\varepsilon|(\Omega)
&\leq J_\varepsilon(\rho_\varepsilon)
     -\varepsilon H(\rho_\varepsilon)\\
&\leq F(\rho_0)-\varepsilon\log|\Omega|
     -\varepsilon H(\rho_\varepsilon)
\leq nD^2.
\end{aligned}
$$
Together with the unit-mass constraint, this gives a
uniform $BV$ bound and hence a uniform $L^r$ bound.
Since $|s\log s|\leq C_r(1+s^r)$ for $s\geq0$,
the entropies $H(\rho_\varepsilon)$ are uniformly bounded.
By compactness, along a sequence $\varepsilon\to0$, $
\rho_\varepsilon\to\rho_* $ strongly in $L^r(\Omega)$, where $\rho_*\in\mathcal A$.
For every $\eta\in\mathcal A$, optimality yields
$$
J(\rho_\varepsilon)
\leq J(\eta)
+\varepsilon\bigl(H(\eta)-H(\rho_\varepsilon)\bigr).
$$
Passing to the limit by lower semicontinuity gives
$J(\rho_*)\leq J(\eta)$, so $\rho_*$ minimizes $J$.

\noindent
\textit{Step 2.} We next establish a positive lower bound for each fixed
$\varepsilon>0$, following the perturbation argument
of~\cite[Proposition~3.1]{CP19}.
For fixed $\varepsilon$, choose $t>0$ such that $
E:=\{\rho_\varepsilon>t\}$, has positive measure and finite relative
perimeter. For $
0<\delta\leq\frac{t|E|}{2|\Omega|}$ and $
m_\delta:=\int_\Omega(\delta-\rho_\varepsilon)_+\,dx$, define $
\eta_\delta
:=\max\{\rho_\varepsilon,\delta\}
-\frac{m_\delta}{|E|}\mathbf1_E$. Indeed, $m_\delta\leq\delta|\Omega|$ implies
$m_\delta/|E|\leq t/2$. Since $\delta\leq t/2$, the sets
$\{\rho_\varepsilon<\delta\}$ and $E=\{\rho_\varepsilon>t\}$
are disjoint. The construction adds mass $m_\delta$
on the former and removes the same amount from the latter.
Moreover,
$$
\frac{m_\delta}{|E|}
\leq\frac{\delta|\Omega|}{|E|}
\leq\frac t2,
$$
so $\eta_\delta\geq t/2$ on $E$.
Truncation preserves $BV$, and $\mathbf1_E\in BV(\Omega)$
because $E$ has finite relative perimeter. Hence
$\eta_\delta\in\mathcal A$
and $\|\eta_\delta-\rho_\varepsilon\|_{L^1(\Omega)}
=2m_\delta$. Since $s\mapsto\max\{s,\delta\}$ is $1$-Lipschitz,
$
\bigl|D\max\{\rho_\varepsilon,\delta\}\bigr|(\Omega)
\leq |D\rho_\varepsilon|(\Omega)$. Thus, by the triangle inequality for the total variation,
$$
\begin{aligned}
|D\eta_\delta|(\Omega)
&\leq
\bigl|D\max\{\rho_\varepsilon,\delta\}\bigr|(\Omega)
+\frac{m_\delta}{|E|}|D\mathbf1_E|(\Omega)\\
&\leq |D\rho_\varepsilon|(\Omega)
+\frac{m_\delta}{|E|}\operatorname{Per}(E;\Omega).
\end{aligned}
$$
By \eqref{eq:transport-L1}, we have $
F(\eta_\delta)-F(\rho_\varepsilon)\leq4nD^2m_\delta$. To estimate the entropy, write $h(s)=s\log s$.
On $E$, set $a:=m_\delta/|E|$. Since
$\eta_\delta=\rho_\varepsilon-a\geq t/2$, we have
$$
\begin{aligned}
h(\eta_\delta)-h(\rho_\varepsilon)
&=-\int_{\eta_\delta}^{\rho_\varepsilon}
  (1+\log s)\,ds\\
&\leq a\bigl[-1-\log(t/2)\bigr]
\leq aC(t),
\end{aligned}
$$
where $C(t):=\max\{0,-1-\log(t/2)\}$.
Integrating over $E$ yields
$
\int_E\bigl[h(\eta_\delta)-h(\rho_\varepsilon)\bigr]\,dx
\leq m_\delta C(t)$. Since $\eta_\delta=\delta$ on
$A_\delta:=\{\rho_\varepsilon<\delta\}$ and
$\eta_\delta=\rho_\varepsilon$ outside $A_\delta\cup E$,
the disjointness of $A_\delta$ and $E$ gives
$$
\begin{aligned}
H(\eta_\delta)-H(\rho_\varepsilon)
&=\int_{A_\delta}
  \bigl[h(\delta)-h(\rho_\varepsilon)\bigr]\,dx
  +\int_E
  \bigl[h(\eta_\delta)-h(\rho_\varepsilon)\bigr]\,dx\\
&\leq(1+\log\delta)
  \int_{A_\delta}(\delta-\rho_\varepsilon)\,dx
  +m_\delta C(t)\\
&=m_\delta\bigl(1+\log\delta+C(t)\bigr).
\end{aligned}
$$ Combining these estimates with optimality gives
\begin{equation}\label{eq:entropy-perturb}
0
\leq J_\varepsilon(\eta_\delta)
     -J_\varepsilon(\rho_\varepsilon)
\leq
m_\delta\left[
4nD^2+\frac{\operatorname{Per}(E;\Omega)}{|E|}
+\varepsilon\bigl(1+\log\delta+C(t)\bigr)
\right].
\end{equation}
For sufficiently small $\delta>0$, the bracket is
strictly negative. Hence $m_\delta=0$, which implies
$\rho_\varepsilon\geq\delta$ almost everywhere and proves
$\operatorname*{ess\,inf}_\Omega\rho_\varepsilon>0$.

\noindent
\textit{Step 3.} To obtain a uniform bound for $b_\varepsilon^-$ along
the selected sequence, choose $0<a<b<c$ such that
$|\{\rho_*>c\}|>0$. Set $v:=\frac12|\{\rho_*>c\}|>0$.
Convergence in measure implies that, for all sufficiently
small $\varepsilon$ and every $t\in[a,b]$, $|\{\rho_\varepsilon>t\}|\geq v>0$
with $v$ independent of $\varepsilon$.
By the coarea formula and the uniform variation bound,
we can choose $t_\varepsilon\in(a,b)$ such that $
E_\varepsilon:=\{\rho_\varepsilon>t_\varepsilon\}$ and $
\operatorname{Per}(E_\varepsilon;\Omega)
\leq\frac{nD^2}{b-a}$.
Thus
$$
\frac{\operatorname{Per}(E_\varepsilon;\Omega)}
     {|E_\varepsilon|}
\leq\frac{nD^2}{(b-a)v}.
$$
Set $c_0:=\min\{1,av/(2|\Omega|)\}$.
Using $E_\varepsilon$ in \eqref{eq:entropy-perturb} and
choosing $
\delta_\varepsilon
:=\min\{c_0,\exp(-C/\varepsilon)\}$
with $C$ sufficiently large and independent of
$\varepsilon$, the bracket is strictly negative.
Therefore, $
\rho_\varepsilon\geq\delta_\varepsilon$ almost everywhere
and hence $
0\leq b_\varepsilon^-
\leq\max\{-\varepsilon\log c_0,C\}$.
This proves the uniform $L^\infty$ bound. Finally, for every $1\leq s<\infty$ and $t\geq0$, $
\bigl(\log\max\{1,t\}\bigr)^s
\leq C_{r,s}(1+t^r)$. The uniform $L^r$ bound therefore yields
$$
\|b_\varepsilon^+\|_{L^s(\Omega)}^s
\leq C_{r,s}\varepsilon^s
\left(|\Omega|+\|\rho_\varepsilon\|_{L^r(\Omega)}^r\right)
\leq C_s\varepsilon^s.
$$
Thus $b_\varepsilon^+\to0$ in every finite $L^s(\Omega)$.
\end{proof}

\subsection{First variation and optimality conditions}

The positivity obtained above allows us to differentiate
the transport terms with respect to the density.
We first establish the required first-variation formula
and then combine it with the entropy derivative
and the total variation inequality to obtain
the regularized optimality conditions.

\begin{lemma}[First variation at a positive density]\label{L5.3}
Let $\sigma\in\A$ have positive essential infimum. For arbitrary
$\mu\in\Prob(K)$,  there exists a unique
$c$-concave potential $f\in C(K)$ satisfying $f(x_*)=0$
and $
W_2^2(\mu,\sigma\,dx)
=\int_K f^c\,d\mu+\int_\Omega f\sigma\,dx$. For any probability density $\eta\in L^1(\Omega)$, set
$\sigma_t=(1-t)\sigma+t\eta$ for $0\leq t<1$. Then
\begin{equation}\label{eq:first-variation}
\lim_{t\to0}
\frac{W_2^2(\mu,\sigma_t\,dx)-W_2^2(\mu,\sigma\,dx)}{t}
=
\int_\Omega f(\eta-\sigma)\,dx.
\end{equation}
Moreover, if $w\in L^\infty(\Omega)$ satisfies
$\int_\Omega w\,dx=0$, then
$$
\lim_{t\to0}
\frac{W_2^2(\mu,(\sigma+tw)\,dx)-W_2^2(\mu,\sigma\,dx)}{t}=
\int_\Omega fw\,dx,
$$
where $\sigma+tw$ is nonnegative for all sufficiently
small $|t|$.
\end{lemma}

\begin{proof}
Let $f$ and $g$ be two $c$-concave optimal potentials associated with $\sigma\,dx$, both normalized by $f(x_*)=g(x_*)=0$.
Since $\sigma\,dx$ is absolutely continuous, Brenier's
theorem gives a unique optimal transport map from
$\sigma\,dx$ to $\mu$. For the cost $|x-y|^2$, this map
is represented by both
$$
T(x)=x-\frac12\nabla f(x)
    =x-\frac12\nabla g(x)
\qquad\text{for }\sigma\,dx\text{-almost every }x.
$$
As $\sigma>0$ almost everywhere, it follows that
$\nabla f=\nabla g$ Lebesgue-almost everywhere in $\Omega$. Poincar\'e's inequality on the connected domain
$\Omega$ implies that $f-g$ is constant almost everywhere.
Continuity extends this equality to $K$,
and the normalization $f(x_*)=g(x_*)=0$
then gives $f=g$ on $K$.

For $0<t<1$, choose a normalized
$c$-concave optimal potential $f_t\in\mathcal K$ for
$(\mu,\sigma_t\,dx)$. By compactness of $\mathcal K$,
every sequence $t_k\to0$ admits a subsequence
along which $f_{t_k}$ converges uniformly.
Since $\sigma_t\to\sigma$ in $L^1(\Omega)$ and the $c$-transform is continuous with respect to
the supremum norm on $C(K)$,
passing to the limit in the dual identity shows that
each such limit is a normalized $c$-concave optimal
potential for $(\mu,\sigma\,dx)$.
By the uniqueness proved above, every limit equals $f$.
Hence $f_t\to f$ uniformly on $K$ as $t\to0$. Using $f$ as a dual competitor at $\sigma_t$ and $f_t$
as a dual competitor at $\sigma$, we obtain
$$
\int_\Omega f(\eta-\sigma)\,dx
\leq
\frac{W_2^2(\mu,\sigma_t\,dx)-W_2^2(\mu,\sigma\,dx)}{t}
\leq
\int_\Omega f_t(\eta-\sigma)\,dx.
$$
Uniform convergence of $f_t$ implies that both bounds
converge to $\int_\Omega f(\eta-\sigma)\,dx$, proving
\eqref{eq:first-variation}.
The same argument applies to $\sigma+tw$ for bounded
zero integral perturbations $w$, since the positive lower
bound of $\sigma$ ensures admissibility for sufficiently
small $|t|$.
\end{proof}

We apply this formula at the regularized minimizer
$\rho_\varepsilon$ and use the representation
of total variation established above.

\begin{proposition}[The approximate Euler-Lagrange equation]
\label{prop:approximate-EL}
Let $f_{i,\varepsilon}\in\K$ be normalized $c$-concave potentials
for $(\mu_i,\rho_\varepsilon\dd x)$ and set
$V_\varepsilon=\sum_i f_{i,\varepsilon}$. There exist
$\lambda_\varepsilon\in\R$ and $z_\varepsilon\in L^\infty$
with $\|z_\varepsilon\|_\infty\leq1$, zero normal trace, and
\begin{align}
 &\diver z_\varepsilon
 =\lambda_\varepsilon-V_\varepsilon
                        -\varepsilon\log\rho_\varepsilon,
 \label{eq:approx-EL}\\
 &\int_\Omega\rho_\varepsilon\diver z_\varepsilon\dd x
 =\TV(\rho_\varepsilon).
 \label{eq:approx-calibration}
\end{align}
The divergence belongs to every finite $L^s(\Omega)$.
\end{proposition}
\begin{proof}
Take $w\in BV(\Omega)\cap L^\infty(\Omega)$ with
$\int_\Omega w\,dx=0$.
Since $\operatorname*{ess\,inf}_\Omega\rho_\varepsilon>0$ (Lemma~\ref{lem:entropy-bounds}), we have $\rho_\varepsilon+tw\in\A$
for all sufficiently small $t>0$.
For such $t$, the triangle inequality gives $
|D(\rho_\varepsilon+tw)|(\Omega)
\leq |D\rho_\varepsilon|(\Omega)+t|Dw|(\Omega)$.
Minimality of $\rho_\varepsilon$ therefore yields
$$
\begin{aligned}
0
&\leq
\frac{J_\varepsilon(\rho_\varepsilon+tw)
      -J_\varepsilon(\rho_\varepsilon)}{t}\\
&\leq
\frac{F(\rho_\varepsilon+tw)-F(\rho_\varepsilon)}{t}
+|Dw|(\Omega)
+\varepsilon
\frac{H(\rho_\varepsilon+tw)-H(\rho_\varepsilon)}{t}.
\end{aligned}
$$
To differentiate the entropy, write $h(s)=s\log s$.
The positive lower bound for $\rho_\varepsilon$
and boundedness of $w$ give, for sufficiently small $t>0$,
$$
\left|
\frac{h(\rho_\varepsilon+tw)-h(\rho_\varepsilon)}{t}
\right|
\leq
C_{\varepsilon,w}
\bigl(1+\log(1+\rho_\varepsilon)\bigr).
$$
The right-hand side is integrable and independent of $t$.
Dominated convergence therefore yields
\[
\lim_{t\downarrow0}
\frac{H(\rho_\varepsilon+tw)-H(\rho_\varepsilon)}{t}
=
\int_\Omega(1+\log\rho_\varepsilon)w\,dx.
\]
Using the first-variation formula for $F$ and
$\int_\Omega w\,dx=0$, we conclude that
\begin{equation}\label{e5.11}
|Dw|(\Omega)
\geq-\int_\Omega
(V_\varepsilon+\varepsilon\log\rho_\varepsilon)w\,dx.
\end{equation}
The constant $1$ in the derivative of $\rho\log\rho$ cancels
because $\int w=0$. Define
$$\lambda_\varepsilon=\frac1{|\Omega|}\int_\Omega(V_\varepsilon+\varepsilon\log\rho_\varepsilon)\dd x,
\qquad
\xi_\varepsilon=\lambda_\varepsilon-V_\varepsilon-\varepsilon\log\rho_\varepsilon.
$$
By construction, $\int_\Omega\xi_\varepsilon\,dx=0$.
For any $w\in BV(\Omega)\cap L^\infty(\Omega)$, apply
the preceding inequality \eqref{e5.11} to $
w-\overline w$ with $\overline w:=\frac1{|\Omega|}\int_\Omega w\,dx$. Subtracting a constant changes neither the total variation
nor the integral against $\xi_\varepsilon$, so $
\int_\Omega\xi_\varepsilon w\,dx\leq|Dw|(\Omega)$. The positive lower bound for $\rho_\varepsilon$ and its
logarithmic integrability imply
$\log\rho_\varepsilon\in L^q(\Omega)$, hence
$\xi_\varepsilon\in L^q(\Omega)$.
Lemma~\ref{lem:TV-polar} therefore yields a vector field
$z_\varepsilon\in L^\infty(\Omega;\mathbb R^d)$ with
$\|z_\varepsilon\|_\infty\leq1$, zero normal trace, and $
\operatorname{div}z_\varepsilon=\xi_\varepsilon$.

It remains to prove \eqref{eq:approx-calibration}.
Let $\rho_0:=|\Omega|^{-1}$ and consider $
\sigma_t:=(1-t)\rho_\varepsilon+t\rho_0$, where $0\leq t<1$. This path lies in $\mathcal A$. Since $\rho_0$ is constant, $|D\sigma_t|(\Omega)
=(1-t)|D\rho_\varepsilon|(\Omega)$.
Lemma~\ref{L5.3} gives
$$
\left.\frac{d}{dt}F(\sigma_t)\right|_{t=0^+}
=\int_\Omega V_\varepsilon(\rho_0-\rho_\varepsilon)\,dx.
$$
For $0<t<1/2$, the positive lower bound for
$\rho_\varepsilon$ and the mean value theorem give
$$
\left|
\frac{h(\sigma_t)-h(\rho_\varepsilon)}{t}
\right|
\leq C_{\varepsilon,r}(1+\rho_\varepsilon^r),
$$
where the right-hand side belongs to $L^1(\Omega)$
and is independent of $t$. Dominated convergence
therefore yields
$$
\left.\frac{d}{dt}H(\sigma_t)\right|_{t=0^+}
=
\int_\Omega(1+\log\rho_\varepsilon)
(\rho_0-\rho_\varepsilon)\,dx.
$$
The constant term vanishes because both densities have
unit mass. By minimality of $\rho_\varepsilon$,
$$
\begin{aligned}
0
&\leq\left.\frac{d}{dt}J_\varepsilon(\sigma_t)\right|_{t=0^+}\\
&=-|D\rho_\varepsilon|(\Omega)
+\int_\Omega
(V_\varepsilon+\varepsilon\log\rho_\varepsilon)
(\rho_0-\rho_\varepsilon)\,dx\\
&=-|D\rho_\varepsilon|(\Omega)
+\int_\Omega\xi_\varepsilon\rho_\varepsilon\,dx.
\end{aligned}
$$
In the last equality, we used
$V_\varepsilon+\varepsilon\log\rho_\varepsilon
=\lambda_\varepsilon-\xi_\varepsilon$,
the equality of masses, and
$\int_\Omega\xi_\varepsilon\rho_0\,dx=0$. Conversely, \eqref{eq:weak-pairing-bound} and
$\operatorname{div}z_\varepsilon=\xi_\varepsilon$ give $
\int_\Omega\xi_\varepsilon\rho_\varepsilon\,dx
\leq|D\rho_\varepsilon|(\Omega)$.
Hence $
\int_\Omega\rho_\varepsilon
\operatorname{div}z_\varepsilon\,dx
=|D\rho_\varepsilon|(\Omega)$, which proves \eqref{eq:approx-calibration}. Finally, \eqref{eq:approx-EL} and the logarithmic integrability estimates imply
$\operatorname{div}z_\varepsilon\in L^s(\Omega)$
for every $1\leq s<\infty$.
\end{proof}

\subsection{The vector-field dual problem}

We now pass to the limit in the regularized
optimality conditions.
This produces a maximizing tuple for the vector-field
dual and yields necessary and sufficient conditions
for a density to minimize $J$.
Let $\mathcal Z$ be the fields $z\in L^\infty(\Omega;\R^d)$ with
$\|z\|_\infty\leq1$, $\diver z\in L^\infty(\Omega)$, and zero
normal trace in the sense of \eqref{eq:weak-zero-trace}.

\begin{theorem}[Vector-field duality and optimality]\label{thm:field-dual}
For arbitrary prescribed measures $\mu_i\in\Prob(K)$,
\begin{equation}\label{eq:field-dual}
 \min_{\rho\in\A}J(\rho)
 =\max\left\{\lambda+\sum_{i=1}^n\int_K f_i^c\dd\mu_i:
 \begin{array}{l}
 f_i\in C(K),\ z\in\mathcal Z,\ \lambda\in\R,\\
 \sum_i f_i+\diver z\geq\lambda\quad\text{a.e.}
 \end{array}\right\}.
\end{equation}
There is a maximizing tuple with $f_i\in\K$. For any fixed
maximizing tuple, set $\beta=\sum_i f_i+\diver z-\lambda$.
A density $\rho\in\A$ is a minimizer if and only if
\begin{align}
 &W_2^2(\mu_i,\rho\,d x)
 =\int_K f_i^c\,d\mu_i+\int_\Omega f_i\rho\,d x\qquad(i=1,\ldots,n),\label{eq:field-OT}\\
&\beta\rho=0\text{ a.e. and } 
 \int_\Omega\rho\diver z\dd x=|D\rho|(\Omega).
 \label{eq:field-complementarity}
\end{align}
In particular, $\beta\in L^\infty$, $\beta\geq0$, and the TV
condition is equivalently $-(z,D\rho)=|D\rho|$ as measures.
\end{theorem}
\begin{proof}
Let $(\boldsymbol f,z,\lambda)$ be any admissible tuple
for the dual problem, and set $V:=\sum_{i=1}^{n}f_i$.
For every $\eta\in\mathcal A$, Kantorovich duality gives
$$
\sum_{i=1}^{n}W_2^2(\mu_i,\eta\,dx)
\geq B(\boldsymbol f)+\int_\Omega V\eta\,dx,
$$
while \eqref{eq:weak-pairing-bound} yields $
|D\eta|(\Omega)
\geq\int_\Omega\eta\,\operatorname{div}z\,dx$. Combining these inequalities with the dual constraint
$V+\operatorname{div}z\geq\lambda$, we obtain
$$
\begin{aligned}
J(\eta)
&\geq B(\boldsymbol f)
+\int_\Omega(V+\operatorname{div}z)\eta\,dx\\
&\geq B(\boldsymbol f)+\lambda\int_\Omega\eta\,dx\\
&=B(\boldsymbol f)+\lambda,
\end{aligned}
$$
where we used $\eta\geq0$ and $\int_\Omega\eta\,dx=1$.
Taking the infimum over $\eta\in\mathcal A$ and the
supremum over all dual admissible tuples shows that the dual supremum does not exceed the primal minimum.

We next pass to the limit in the regularized optimality
conditions. Along the sequence provided by
Lemma~\ref{lem:entropy-bounds},
$\rho_\varepsilon\to\rho_*$ strongly in $L^r(\Omega)$,
where $\rho_*$ minimizes $J$.
Since $f_{i,\varepsilon}\in\mathcal K$ and $\mathcal K$
is compact in $C(K)$, we may extract a subsequence such that $
f_{i,\varepsilon}\rightarrow f_i\in\mathcal K$ uniformly on $K$, for $i=1,\ldots,n$. Consequently, $V_\varepsilon:=\sum_i f_{i,\varepsilon}$
converges uniformly to $V:=\sum_i f_i$. Using $\varepsilon\log\rho_\varepsilon
=b_\varepsilon^+-b_\varepsilon^-$, we can write
\[
\lambda_\varepsilon
=\frac1{|\Omega|}
\int_\Omega
\bigl(V_\varepsilon+b_\varepsilon^+-b_\varepsilon^-\bigr)\,dx.
\]
The uniform bound for $V_\varepsilon$ and
\eqref{eq:entropy-log-bounds} therefore imply that
$(\lambda_\varepsilon)$ is bounded.
Since $\|z_\varepsilon\|_\infty\leq1$ and
$(b_\varepsilon^-)$ is bounded in $L^\infty(\Omega)$,
a further subsequence satisfies
$
\lambda_\varepsilon\to\lambda$, $z_\varepsilon\rightharpoonup^*z$ in $L^\infty(\Omega;\mathbb R^d)$, and
$b_\varepsilon^-\rightharpoonup^*\beta$ in $L^\infty(\Omega)$. These limits satisfy $\|z\|_\infty\leq1$ and $\beta\geq0$
almost everywhere. By \eqref{eq:approx-EL},
$$
\xi_\varepsilon:=\operatorname{div}z_\varepsilon
=\lambda_\varepsilon-V_\varepsilon
-b_\varepsilon^++b_\varepsilon^-.
$$
Since $b_\varepsilon^+\to0$ strongly in $L^q(\Omega)$,
the preceding convergences give
$$
\xi_\varepsilon\rightharpoonup
\lambda-V+\beta=:\xi
\qquad\text{in }L^q(\Omega).
$$
For every $\varphi\in C^1(\overline\Omega)$,
the zero normal trace condition gives
$
\int_\Omega\xi_\varepsilon\varphi\,dx
=-\int_\Omega z_\varepsilon\cdot\nabla\varphi\,dx$.
Passing to the limit yields
$
\int_\Omega\xi\varphi\,dx
=-\int_\Omega z\cdot\nabla\varphi\,dx$.
Thus $\operatorname{div}z=\xi$ in the distributional sense,
and $z$ has zero normal trace.
Moreover, $\xi=\lambda-V+\beta\in L^\infty(\Omega)$,
so $z\in\mathcal Z$ and $
V+\operatorname{div}z=\lambda+\beta\geq\lambda$ almost everywhere in $\Omega$. Hence $(\boldsymbol f,z,\lambda)$ is admissible for
the dual problem. The elementary bound $-t\log t\leq e^{-1}$ for $0<t\leq1$ gives
$$
 0\leq\int_\Omega\rho_\varepsilon b_\varepsilon^-\dd x\leq\frac{\varepsilon|\Omega|}{e}.
$$
Strong $L^1$ convergence of $\rho_\varepsilon$ and weak-* convergence
of $b_\varepsilon^-$ imply $\int\rho_*\beta=0$, hence
$\rho_*\beta=0$.
Moreover, strong convergence in $L^r$ and weak convergence in $L^q$
yield
$$
\TV(\rho_*)\leq\liminf_\varepsilon\TV(\rho_\varepsilon)=\lim_\varepsilon\int_\Omega\rho_\varepsilon\xi_\varepsilon\dd x
 =\int_\Omega\rho_*\diver z\dd x\leq\TV(\rho_*).
$$
All terms are therefore equal.
For each $i$, the estimate used to prove the
$L^1$ continuity of the transport costs gives
$$
\left|
W_2^2(\mu_i,\rho_\varepsilon\,dx)
-W_2^2(\mu_i,\rho_*\,dx)
\right|
\leq2D^2\|\rho_\varepsilon-\rho_*\|_{L^1(\Omega)}
\rightarrow0.
$$
For each $i$, the transport estimate gives $
W_2^2(\mu_i,\rho_\varepsilon\,dx)\rightarrow
W_2^2(\mu_i,\rho_*\,dx)$.
Moreover, uniform convergence of $f_{i,\varepsilon}$
and stability of the $c$-transform imply $
\int_K f_{i,\varepsilon}^c\,d\mu_i
\rightarrow
\int_K f_i^c\,d\mu_i$.
Since $\rho_\varepsilon$ has unit mass,
$$
\begin{aligned}
\left|
\int_\Omega f_{i,\varepsilon}\rho_\varepsilon\,dx
-\int_\Omega f_i\rho_*\,dx
\right|
&\leq
\|f_{i,\varepsilon}-f_i\|_\infty
+\|f_i\|_\infty
 \|\rho_\varepsilon-\rho_*\|_{L^1(\Omega)}
\rightarrow0.
\end{aligned}
$$
Passing to the limit in the optimality identity
\[
W_2^2(\mu_i,\rho_\varepsilon\,dx)
=
\int_K f_{i,\varepsilon}^c\,d\mu_i
+\int_\Omega f_{i,\varepsilon}\rho_\varepsilon\,dx
\]
therefore yields
\[
W_2^2(\mu_i,\rho_*\,dx)
=
\int_K f_i^c\,d\mu_i
+\int_\Omega f_i\rho_*\,dx.
\]
Thus $(f_i^c,f_i)$ is optimal for
$(\mu_i,\rho_*\,dx)$. Multiplying
$V+\diver z=\lambda+\beta$ by $\rho_*$ now gives
$J(\rho_*)=B(\boldsymbol f)+\lambda$.
This proves strong duality and attainment.

We now characterize all primal minimizers.
Fix any dual maximizer $(\boldsymbol f,z,\lambda)$ and set
$$
V:=\sum_{i=1}^{n}f_i,
\qquad
\beta:=V+\operatorname{div}z-\lambda\geq0.
$$
For every $\rho\in\mathcal A$, the unit-mass constraint gives
$$
\begin{aligned}
J(\rho)-\lambda-B(\boldsymbol f)
&=\sum_{i=1}^{n}
\left[
W_2^2(\mu_i,\rho\,dx)
-\int_K f_i^c\,d\mu_i
-\int_\Omega f_i\rho\,dx
\right]\\
&\quad+
\left[
|D\rho|(\Omega)
-\int_\Omega\rho\,\operatorname{div}z\,dx
\right]
+\int_\Omega\beta\rho\,dx.
\end{aligned}
$$
The transport terms are nonnegative by Kantorovich
duality, the second bracket is nonnegative by
\eqref{eq:weak-pairing-bound}, and the last integral is
nonnegative because $\beta,\rho\geq0$. Since $\lambda+B(\boldsymbol f)=\min_{\mathcal A}J$,
a density $\rho$ is a minimizer if and only if every
term on the right vanishes. This is equivalent to
\eqref{eq:field-OT} together with
$$
\int_\Omega\rho\,\operatorname{div}z\,dx
=|D\rho|(\Omega),
\qquad
\beta\rho=0
\quad\text{almost everywhere}.
$$
Thus these conditions are necessary and sufficient
for every primal minimizer. Finally, the estimate $|(z,D\rho)|\leq|D\rho|$ shows
that $|D\rho|+(z,D\rho)$ is a nonnegative measure. By \eqref{eq:weak-pairing-bound} and the optimality condition,
$$
|D\rho|(\Omega)+\int_\Omega d(z,D\rho)
=
|D\rho|(\Omega)
-\int_\Omega\rho\,\operatorname{div}z\,dx
=0.
$$
Hence $
-(z,D\rho)=|D\rho|$
as measures. Conversely, this measure identity implies
the integral equality by the same integration-by-parts
formula.
\end{proof}

\begin{remark}[Relation between the two dual formulations]
Let $(\boldsymbol f,z,\lambda)$ satisfy the constraints
in \eqref{eq:field-dual}, and write
$V=\sum_i f_i$ and $\beta=V+\operatorname{div}z-\lambda\geq0$.
For every set $E\subset\Omega$ of finite relative perimeter,
applying \eqref{eq:weak-pairing-bound} to $\mathbf1_E$ gives $
\operatorname{Per}(E;\Omega)
\geq\int_E\operatorname{div}z\,dx$.
Consequently,
$$
\operatorname{Per}(E;\Omega)+\int_E(V-\lambda)\,dx
\geq\int_E\beta\,dx\geq0.
$$
Thus $(\boldsymbol f,\lambda)$ also satisfies the constraints
in \eqref{eq:dual-sets}. Both dual problems attain the same
value, namely $\min_{\rho\in\mathcal A}J(\rho)$.
\end{remark}

\section{Bounds, one-dimensional reduction and stability}
\label{sec:qualitative}
We next study further properties of the minimizers
and of the variational problem.
We prove an upper bound for the minimizing density,
identify an exact TV-JKO formulation in one dimension,
and establish stability under approximation
of the prescribed measures.

\subsection{Maximum principle}

In this subsection, assume additionally that
$\Omega$ is convex.
We show that the minimizing density is bounded above
by the largest essential supremum of the prescribed
densities.
The proof combines the projection argument
of~\cite[Section~4]{CP19} with the BV estimate
of~\cite{DMSV16}.

\begin{theorem}[Maximum principle]\label{thm:maximum}
Assume that $\Omega$ is convex and that
$\mu_i=\rho_i\,dx$ with $\rho_i\in L^\infty(\Omega)$
for every $i=1,\ldots,n$. Let $\rho_*$ be the unique
minimizer of $J(\rho)$.
Then
\begin{equation}\label{eq:maximum-principle}
\|\rho_*\|_{L^\infty(\Omega)}
\leq
\max_{1\leq i\leq n}\|\rho_i\|_{L^\infty(\Omega)}.
\end{equation}
\end{theorem}

\begin{proof}
Set $M=\max_i\|\rho_i\|_\infty$ and let
$$\mathcal K_M=\{\eta\,d x:\eta\geq0,\ \int_\Omega\eta\,d x=1,\ \eta\leq M\text{ a.e.}\}.
$$
Since each $\rho_i$ has unit mass, $M\geq|\Omega|^{-1}$, so the uniform probability
measure belongs to $\mathcal K_M$.
Moreover, $\mathcal K_M$ is the set of probability
measures on $K$ dominated by
$M\mathbf1_\Omega\,dx$.
This domination is preserved under narrow convergence.
Thus $\mathcal K_M$ is closed in the compact space
$\mathcal P(K)$, and hence narrowly compact. Let $\widehat\rho\,dx$ be a probability measure in
$\mathcal K_M$ minimizing the Wasserstein distance to $\rho_*\,dx$, that is, $
\widehat\rho\,dx
\in\operatorname*{arg\,min}_{\nu\in\mathcal K_M}
W_2^2(\rho_*\,dx,\nu)$. Define the convex lower semicontinuous function
$$
h_M(s):=
\begin{cases}
0, & 0\leq s\leq M,\\
+\infty, & s>M.
\end{cases}
$$
The projection $\widehat\rho$ minimizes
$\frac12W_2^2(\rho_*\,dx,\eta\,dx)+\int_\Omega h_M(\eta)\,dx
$ over probability densities $\eta$.
Since $\Omega$ is convex and $\rho_*\in BV(\Omega)$,
\cite[Theorem~1.1]{DMSV16} applies and gives
\begin{equation}\label{eq:projection-TV}
 |D\widehat\rho|(\Omega)\leq|D\rho_*|(\Omega).
\end{equation}
For completeness, we recall the transport comparison used in
\cite[Lemma 4.1]{CP19}. Fix $i$, and let $T_0,T_1$ be the optimal
maps from $\rho_*\,d x$ to
$\mu_i$ and $\widehat\rho\,d x$, respectively. Define the generalized
geodesic
$$
 \nu_t=((1-t)T_0+tT_1)_\#(\rho_*\,d x),\qquad 0\leq t\leq1.
$$
For every $t\in[0,1]$, the base measure $\rho_*\,dx$ is absolutely continuous,
and both endpoint densities are bounded above by $M$.
The $L^\infty$ bound along generalized geodesics
\cite[Section~4.1]{CP19} therefore yields
$$
\nu_t=\eta_t\,dx,
\qquad
0\leq\eta_t\leq M
\quad\text{a.e.},\qquad t\in[0,1].
$$
Convexity of $\Omega$ ensures that the interpolated maps
take values in $\Omega$ almost everywhere.
Consequently, $\nu_t\in K_M$.
Using the transport plan
$(\mathrm{id},(1-t)T_0+tT_1)_\#(\rho_*\,dx)$, we obtain
\begin{align*}
W_2^2(\rho_*\,dx,\nu_t)
&\leq \int_\Omega
\bigl|x-(1-t)T_0(x)-tT_1(x)\bigr|^2\rho_*(x)\,dx\\
&=(1-t)W_2^2(\rho_*\,dx,\mu_i)
+tW_2^2(\rho_*\,dx,\widehat\rho\,dx)\\
&\quad-t(1-t)\int_\Omega|T_0-T_1|^2\rho_*\,dx\\
&\leq(1-t)W_2^2(\rho_*\,dx,\mu_i)
+tW_2^2(\rho_*\,dx,\widehat\rho\,dx)\\
&\quad-t(1-t)W_2^2(\mu_i,\widehat\rho\,dx),
\end{align*}
where the last inequality follows from
$(T_0,T_1)_\#(\rho_*\,dx)\in\Pi(\mu_i,\widehat\rho\,dx)$.
Projection optimality bounds the left-hand side below by
$W_2^2(\rho_*\,d x,\widehat\rho\,d x)$. Dividing the resulting
inequality by $1-t$ and letting $t\to1$ yields
$$
 W_2^2(\mu_i,\widehat\rho\,d x)
 \leq W_2^2(\mu_i,\rho_*\,d x)-W_2^2(\rho_*\,d x,\widehat\rho\,d x).
$$
Summing over $i$ and using \eqref{eq:projection-TV} gives
$$
 J(\widehat\rho)\leq J(\rho_*)-nW_2^2(\rho_*\,d x,\widehat\rho\,d x).
$$
Optimality of $\rho_*$ forces $\rho_*=\widehat\rho$ a.e., proving
\eqref{eq:maximum-principle}.
\end{proof}

\begin{remark}
The convexity assumption and the quadratic cost are used in this
projection argument. We do not infer a higher-dimensional minimum
principle, or a corresponding maximum principle for every power cost,
from this proof.
\end{remark}

\subsection{TV-JKO formulation in one dimension}

In one dimension, the representation of Wasserstein
distances through generalized inverses of distribution
functions yields an exact variance identity.
Consequently, the sum of transport costs can be
replaced, up to an additive constant, by a single
transport cost to the ordinary Wasserstein barycenter.

\begin{proposition}\label{prop:1D}
Let $\Omega=(a,b)\subset\mathbb R$ and
$\mu_1,\ldots,\mu_n\in\mathcal P([a,b])$. For each $i$, let $Q_i:(0, 1)\to[a, b]$ denote the generalized inverse
of the distribution function of $\mu_i$. Define
$$
 \overline Q=\frac1n\sum_{i=1}^nQ_i,\qquad
 \overline\mu=(\overline Q)_\#\mathcal L^1_{(0,1)},
$$
where $\mathcal L^1_{(0,1)}$ denotes Lebesgue probability measure on $(0,1)$.
Then, for any $\nu\in\Prob([a,b])$, $\overline\mu$ minimizes
$\nu\mapsto\sum_{i=1}^{n}W_2^2(\mu_i,\nu)$
over $\mathcal P([a,b])$, and ,
\begin{equation}\label{eq:1D-variance}
 \sum_iW_2^2(\mu_i,\nu)
=\sum_iW_2^2(\mu_i,\overline\mu)+nW_2^2(\overline\mu,\nu).
\end{equation}
Consequently the TV-regularized problem is exactly, up to an additive
constant, the problem with the single reference measure $\overline\mu$
\begin{equation}\label{eq:1D-JKO}
 \min_{\rho\in\A}\{nW_2^2(\overline\mu,\rho\dd x)+|D\rho|(\Omega)\}.
\end{equation}
Its TV-JKO time step is $\tau=1/(2n)$.
\end{proposition}
\begin{proof}
Since each $Q_i$ is nondecreasing, so is $\overline Q=n^{-1}\sum_{i=1}^n Q_i$.
By the definition of $\overline\mu$, the generalized inverse $Q_{\overline\mu}$ of the
distribution function of $\overline\mu$ satisfies $
Q_{\overline\mu}=\overline Q$ a.e.~on $(0,1)$. For any $\nu\in\mathcal P([a,b])$, let $Q_\nu$ denote
the generalized inverse of its distribution function.
In one dimension, the coupling
$(Q_i,Q_\nu)_\#\mathcal L^1_{(0,1)}$ is optimal
for the quadratic cost, since both maps are nondecreasing.
Consequently,
$$
W_2^2(\mu_i,\nu)
=\int_0^1|Q_i(s)-Q_\nu(s)|^2\,ds.
$$
Writing $Q_i-Q_\nu=(Q_i-\overline Q)+(\overline Q-Q_\nu)$
and using $\sum_{i=1}^n(Q_i-\overline Q)=0$, we obtain
\begin{align*}
\sum_{i=1}^nW_2^2(\mu_i,\nu)
&=\sum_{i=1}^n\int_0^1|Q_i-\overline Q|^2\,ds
+n\int_0^1|\overline Q-Q_\nu|^2\,ds\\
&=\sum_{i=1}^nW_2^2(\mu_i,\overline\mu)
+nW_2^2(\overline\mu,\nu).
\end{align*}
The last term is nonnegative and vanishes precisely
when $\nu=\overline\mu$. Hence $\overline\mu$ is the
unique minimizer of $\nu\mapsto\sum_iW_2^2(\mu_i,\nu)$.
Taking $\nu=\rho\,dx$ and adding $|D\rho|(\Omega)$ yields
$$
J(\rho)
=C+nW_2^2(\overline\mu,\rho\,dx)+|D\rho|(\Omega),
\qquad
C:=\sum_{i=1}^nW_2^2(\mu_i,\overline\mu).
$$
Since $C$ is independent of $\rho$, the two minimization
problems have the same minimizers, and their minimum
values differ by $C$.
Finally, comparison with the TV-JKO functional $\frac{1}{2\tau}W_2^2(\overline\mu,\rho\,dx)
+|D\rho|(\Omega)
$ gives $\tau=1/(2n)$.
\end{proof}

\begin{remark}
The identity \eqref{eq:1D-variance} follows from the
one-dimensional formula expressing $W_2^2(\mu,\nu)$
as the squared $L^2(0,1)$ distance between the generalized
inverses of their distribution functions. In higher dimensions, optimal couplings generally do not
admit such a common representation, so this identity
and the resulting reduction to a TV-JKO step do not
hold in general.
\end{remark}

\subsection{Stability}

We examine how the minimum value and the minimizing
densities depend on the prescribed measures.
The estimate below controls the change in the
minimum value, while BV compactness yields
convergence of minimizers.
Write $J_{\boldsymbol\mu}$ for \eqref{eq:primal} with data
$\boldsymbol\mu=(\mu_1,\ldots,\mu_n)$, and
$\mathsf M(\boldsymbol\mu)=\min_\A J_{\boldsymbol\mu}$.
Let $W_1(\mu,\nu)=\min_{\gamma\in\Pi(\mu,\nu)}
\int|x-y|\dd\gamma$.

\begin{proposition}\label{prop:stability}
For any $\boldsymbol\mu,\boldsymbol\nu\in\mathcal P(K)^n$,
\begin{equation}\label{eq:value-stability}
|\mathsf M(\boldsymbol\mu)-\mathsf M(\boldsymbol\nu)|
\leq 2D\sum_{i=1}^n W_1(\mu_i,\nu_i).
\end{equation}
Suppose that $\mu_i^k$ converges narrowly to $\mu_i$
for each $i$, and let $\rho_k$ minimize
$J_{\boldsymbol\mu^k}$.
Then $(\rho_k)$ is relatively compact in $L^1(\Omega)$,
and every limit $\rho$ minimizes $J_{\boldsymbol\mu}$.
Moreover, along any subsequence $\rho_{k_j}\to\rho$
in $L^1(\Omega)$,
\begin{equation}\label{eq:strict-BV-stability}
|D\rho_{k_j}|(\Omega)\rightarrow|D\rho|(\Omega).
\end{equation}
If at least one limiting measure $\mu_i$ is absolutely
continuous, the limiting problem has a unique minimizer,
and both convergences hold for the whole sequence.
\end{proposition}
\begin{proof}
Fix $\eta\in\mathcal A$. For each $i$, choose an optimal
Kantorovich potential $f_i$ such that
$$
W_2^2(\mu_i,\eta\,dx)
=\int_K f_i^c\,d\mu_i+\int_\Omega f_i\eta\,dx.
$$
Using $f_i$ in the dual problem for $(\nu_i,\eta\,dx)$
and choosing a $W_1$-optimal coupling
$\gamma_i\in\Pi(\mu_i,\nu_i)$, we obtain
\begin{align*}
W_2^2(\mu_i,\eta\,dx)-W_2^2(\nu_i,\eta\,dx)
&\leq \int_K f_i^c\,d(\mu_i-\nu_i)\\
&=\int_{K\times K}\bigl(f_i^c(x)-f_i^c(y)\bigr)
\,d\gamma_i(x,y)\\
&\leq 2D\,W_1(\mu_i,\nu_i),
\end{align*}
since $f_i^c$ is $2D$-Lipschitz.
Interchanging $\mu_i$ and $\nu_i$ gives the same bound
for the absolute difference. Summing over $i$ and
noting that the total variation terms cancel yields
$$
|J_{\boldsymbol\mu}(\eta)-J_{\boldsymbol\nu}(\eta)|
\leq 2D\sum_{i=1}^nW_1(\mu_i,\nu_i).
$$
The bound is independent of $\eta$. Taking infima
over $\eta\in\mathcal A$ in the two corresponding
one-sided inequalities proves \eqref{eq:value-stability}. Let $\rho_0=|\Omega|^{-1}$. By the minimality of $\rho_k$,
$$
|D\rho_k|(\Omega)
\leq J_{\boldsymbol\mu^k}(\rho_k)
\leq J_{\boldsymbol\mu^k}(\rho_0)
\leq nD^2.
$$
Together with $\|\rho_k\|_{L^1(\Omega)}=1$, this gives
a uniform bound in $BV(\Omega)$. Hence, after passing
to a subsequence, $\rho_k\to\rho$ in $L^1(\Omega)$.
Nonnegativity and unit mass pass to the limit, and
lower semicontinuity gives $\rho\in\mathcal A$. Since $K$ is compact, the narrow convergences
$\mu_i^k\to\mu_i$ and $\rho_k\,dx\to\rho\,dx$
also hold in $W_2$. The triangle inequality therefore gives
$$
\left|W_2(\mu_i^k,\rho_k\,dx)
      -W_2(\mu_i,\rho\,dx)\right|
\leq W_2(\mu_i^k,\mu_i)
    +W_2(\rho_k\,dx,\rho\,dx)
\rightarrow0.
$$
Moreover, \eqref{eq:value-stability} implies
$\mathsf M(\boldsymbol\mu^k)\to\mathsf M(\boldsymbol\mu)$.
Combining these facts with lower semicontinuity of
the total variation, we obtain
$$
\mathsf M(\boldsymbol\mu)
\leq J_{\boldsymbol\mu}(\rho)
\leq\liminf_{k\to\infty}
J_{\boldsymbol\mu^k}(\rho_k)
=\mathsf M(\boldsymbol\mu).
$$
Thus $\rho$ minimizes $J_{\boldsymbol\mu}$. Consequently,
$$
|D\rho_k|(\Omega)
=\mathsf M(\boldsymbol\mu^k)-\sum_{i=1}^nW_2^2(\mu_i^k,\rho_k\,dx)\rightarrow
\mathsf M(\boldsymbol\mu)-\sum_{i=1}^nW_2^2(\mu_i,\rho\,dx)
=|D\rho|(\Omega).
$$
The same argument applies to any $L^1$-convergent
subsequence. If the limiting minimizer is unique,
relative compactness implies that the whole sequence
converges to it in $L^1(\Omega)$, and the preceding
identity gives convergence of the total variations.
\end{proof}

\section{Extension to $p$-Wasserstein costs for $p>1$}\label{sec:extensions}
Brizzi and Portinale~\cite{BP26} establish $L^q$ estimates
for unregularized $p$-Wasserstein barycenters under
geometric assumptions on the supports, and show that
such integrability need not be preserved in general
when $p\neq2$.
Here we consider the problem with a total variation
penalty and extend the existence, uniqueness, and
perimeter duality results to the cost
$c_p(x,y)=|x-y|^p$, where $p>1$.
The corresponding optimality conditions also yield
a transport-map representation when a prescribed
measure is absolutely continuous.
We state these conclusions together and explain
the changes required in the quadratic-cost arguments. 

For $p>1$, consider
$$
J_p(\rho)
=
\sum_{i=1}^nW_p^p(\mu_i,\rho\,dx)+|D\rho|(\Omega),
\qquad \rho\in\A.
$$
Recall that $\Prob_p(K)=\Prob(K)$ because $K$ is compact.

\begin{theorem}\label{thm:p}
The problem $\min_{\rho\in\A}J_p(\rho)$ has a solution for arbitrary
$\mu_i\in\Prob_p(K)$. The minimizer is unique if at least one $\mu_i$
is absolutely continuous with respect to Lebesgue measure. It also has the attained dual representation
\begin{equation}\label{eq:pdual}
 \min_{\rho\in\A}J_p(\rho)
 =\max_{\boldsymbol f\in\K_p^n}
   \left\{\sum_{i=1}^n\int_K f_i^{c_p}\dd\mu_i
            +m\left(\sum_i f_i\right)\right\},
\end{equation}
where
$\K_p=\{f\in C(K):f(x_*)=0,\ \Lip(f)\leq pD^{p-1}\}$. 
For any dual maximizer $\boldsymbol f\in\mathcal K_p^n$, set $
V:=\sum_{i=1}^n f_i$, $\lambda:=m(V)$. The optimality characterization and superlevel-set
identities established for the quadratic cost hold
with this value of $\lambda$, with $W_2^2$ replaced
by $W_p^p$ and $c$ replaced by $c_p$. If $\mu_j$ is absolutely continuous for some
$j\in\{1,\ldots,n\}$, let $\rho$ denote the unique
minimizer of $J_p$.
For every dual maximizer $\boldsymbol f\in\mathcal K_p^n$,
\begin{equation}\label{eq:p-map}
\rho\,dx=(T_j)_\#\mu_j,
\qquad
T_j(x)=
x-(\nabla h_p)^{-1}
       \bigl(\nabla f_j^{c_p}(x)\bigr)
\quad\text{for }\mu_j\text{-a.e. }x,
\end{equation}
where $h_p(v)=|v|^p$.
Here $(\nabla h_p)^{-1}(0)=0$, and
$$
(\nabla h_p)^{-1}(w)
=
p^{-1/(p-1)}|w|^{(2-p)/(p-1)}w,
\qquad w\in\mathbb R^d\setminus\{0\}.
$$
\end{theorem}
\begin{proof}
The constant competitor gives $J_p(\rho_0)\leq nD^p$. Therefore every bounded-energy minimizing sequence has a subsequence converging in
$L^1$ to an element of $\A$. On the compact space $K$, this implies
convergence in $W_p$. Lower semicontinuity of TV proves existence. Suppose that $\mu_j$ is absolutely continuous.
Since $h_p(v)=|v|^p$ is differentiable and strictly convex
for $p>1$, the optimal plan from $\mu_j$ to any
$\nu\in\mathcal P(K)$ is unique and induced by a map;
see~\cite[Theorem~1.17 and Remark~1.18]{San15}.
The argument in Theorem~\ref{thm:existence}, based on
convex combinations of optimal plans, therefore shows
that $\nu\mapsto W_p^p(\mu_j,\nu)$ is strictly convex
under linear interpolation of measures.
The remaining transport terms and the total variation
are convex, so $J_p$ has a unique minimizer.

For $x,y,y'\in K$,
$$
\bigl||x-y|^p-|x-y'|^p\bigr|
\leq pD^{p-1}|y-y'|.
$$
Thus every $c_p$-transform is $pD^{p-1}$-Lipschitz.
Replacing a potential $f$ by $f^{c_pc_p}$ does not
decrease the Kantorovich objective, and subtracting
$f^{c_pc_p}(x_*)$ leaves it unchanged.
Hence the dual potentials may be restricted to
$\mathcal K_p$, which is compact in $C(K)$ and convex. Set
$$
L_p(\rho,\boldsymbol f)
=
|D\rho|(\Omega)
+\int_\Omega\rho\sum_{i=1}^n f_i\,dx
+\sum_{i=1}^n\int_K f_i^{c_p}\,d\mu_i.
$$
Kantorovich duality yields $
J_p(\rho)=\max_{\boldsymbol f\in\mathcal K_p^n}
L_p(\rho,\boldsymbol f)$. For fixed $\boldsymbol f$, $L_p$ is convex and
lower semicontinuous in $\rho$ with respect to
$L^1(\Omega)$ convergence. For fixed $\rho$, it is
concave and continuous in $\boldsymbol f$, using $
\|f^{c_p}-g^{c_p}\|_\infty\leq\|f-g\|_\infty$. Sion's minimax theorem therefore gives
$$
\begin{aligned}
\min_{\rho\in\mathcal A}J_p(\rho)
&=\sup_{\boldsymbol f\in\mathcal K_p^n}
  \inf_{\rho\in\mathcal A}L_p(\rho,\boldsymbol f)\\
&=\sup_{\boldsymbol f\in\mathcal K_p^n}
\left\{
\sum_{i=1}^n\int_K f_i^{c_p}\,d\mu_i
+m\left(\sum_{i=1}^n f_i\right)
\right\}.
\end{aligned}
$$
The objective is continuous by the preceding estimate
and the Lipschitz continuity of $m$.
Compactness of $\mathcal K_p^n$ ensures that the
supremum is attained, proving \eqref{eq:pdual}.

Let $\rho$ be a primal minimizer and $\boldsymbol f$
a dual maximizer, and set $V=\sum_{i=1}^n f_i$.
The duality identity gives
$$
\begin{aligned}
0={}&
\sum_{i=1}^n
\left[
W_p^p(\mu_i,\rho\,dx)
-\int_K f_i^{c_p}\,d\mu_i
-\int_\Omega f_i\rho\,dx
\right]\\
&+|D\rho|(\Omega)+\int_\Omega V\rho\,dx-m(V).
\end{aligned}
$$
Each term is nonnegative by Kantorovich duality
and the definition of $m(V)$, so every term vanishes. Conversely, these equalities make the primal objective
equal to the attained dual value, and hence imply
optimality of $\rho$.
This gives the transport optimality identities and
$$
|D\rho|(\Omega)+\int_\Omega V\rho\,dx=m(V).
$$
Writing $E_t=\{\rho>t\}$, the coarea and layer-cake
formulas imply
\[
\int_0^\infty
\left[
\operatorname{Per}(E_t;\Omega)
+\int_{E_t}V\,dx-m(V)|E_t|
\right]dt=0.
\]
The integrand is nonnegative by
Lemma~\ref{lem:coarea-dual}. Hence, for almost every
$t>0$ with $|E_t|>0$,
$$
\frac{\operatorname{Per}(E_t;\Omega)+\int_{E_t}V\,dx}
{|E_t|}
=m(V).
$$

Finally, suppose that $\mu_j$ is absolutely continuous,
and let $T_j$ be the optimal transport map from
$\mu_j$ to $\rho\,dx$. The transport optimality identity
implies
$$
f_j^{c_p}(x)+f_j(T_j(x))
=h_p(x-T_j(x))
\quad\text{for }\mu_j\text{-a.e. }x.
$$
At any such interior point where $f_j^{c_p}$ is
differentiable, the function
$$
x'\longmapsto h_p(x'-T_j(x))-f_j(T_j(x))-f_j^{c_p}(x')
$$
is nonnegative and vanishes at $x'=x$. Therefore,
$$
\nabla f_j^{c_p}(x)=\nabla h_p(x-T_j(x)).
$$
Since $f_j^{c_p}$ is Lipschitz and $\mu_j$ is absolutely
continuous, this holds for $\mu_j$-almost every $x$.
Applying $(\nabla h_p)^{-1}$ proves \eqref{eq:p-map}.
\end{proof}

\begin{remark}[The case $p=1$]
Existence and the attained dual remain valid for $p=1$, using the
normalized $1$-Lipschitz class of potentials. The strict-convexity
argument and the inverse-gradient formula above do not extend to this
case. Their failure alone does not imply nonuniqueness of the
TV-regularized barycenter.
\end{remark}


\section*{Acknowledgements}
This work was supported by the National Natural Science Foundation of
China (Grant Nos. 12301661 and 12431018).



\begin{thebibliography}{99}
\small
\setlength{\itemsep}{2pt}

\bibitem{bib1}
M.~Agueh and G.~Carlier,
Barycenters in the Wasserstein space,
\emph{SIAM J. Math. Anal.} \textbf{43} (2011), 904--924.

\bibitem{AFP00}
L.~Ambrosio, N.~Fusco and D.~Pallara,
\emph{Functions of Bounded Variation and Free Discontinuity Problems},
Oxford University Press, Oxford, 2000.

\bibitem{bib4}
L.~Ambrosio, N.~Gigli and G.~Savar\'e,
\emph{Gradient Flows in Metric Spaces and in the Space of Probability Measures},
Birkh\"auser, Basel, 2005.

\bibitem{Anz83}
G.~Anzellotti,
Pairings between measures and bounded functions and compensated compactness,
\emph{Ann. Mat. Pura Appl.} \textbf{135} (1983), 293--318.


\bibitem{bib6}
J.~Bigot, E.~Cazelles and N.~Papadakis,
Penalization of barycenters in the Wasserstein space,
\emph{SIAM J. Math. Anal.} \textbf{51} (2019), 2261--2285.

\bibitem{bib10}
N.~Bonneel, J.~Rabin, G.~Peyr\'e and H.~Pfister,
Sliced and Radon Wasserstein barycenters of measures,
\emph{J. Math. Imaging Vision} \textbf{51} (2015), 22--45.

\bibitem{BP26}
C.~Brizzi and L.~Portinale,
On the $q$-integrability of $p$-Wasserstein barycenters, preprint
\href{https://arxiv.org/abs/2602.18293}{arXiv:2602.18293}, 2026.

\bibitem{bib19}
G.~Carlier, K.~Eichinger and A.~Kroshnin,
Entropic-Wasserstein barycenters: PDE characterization, regularity, and CLT,
\emph{SIAM J. Math. Anal.} \textbf{53} (2021), 5880--5914.

\bibitem{CP19}
G.~Carlier and C.~Poon,
On the total variation Wasserstein gradient flow and the TV-JKO scheme,
\emph{ESAIM Control Optim. Calc. Var.} \textbf{25} (2019), 42.

\bibitem{Chizat26}
L.~Chizat,
Doubly regularized entropic Wasserstein barycenter,
\emph{Found. Comput. Math.} \textbf{26} (2026), 2091--2128.

\bibitem{DMSV16}
G.~De Philippis, A.~R.~M\'esz\'aros, F.~Santambrogio and B.~Velichkov,
BV estimates in optimal transportation and applications,
\emph{Arch. Ration. Mech. Anal.} \textbf{219} (2016), 829--860.


\bibitem{bib13}
A.~Figalli and F.~Glaudo,
\emph{An Invitation to Optimal Transport, Wasserstein Distances, and Gradient Flows},
European Mathematical Society, Berlin, 2021.

\bibitem{bib5}
W.~Gangbo and A.~\'Swi\c{e}ch,
Optimal maps for the multidimensional Monge-Kantorovich problem,
\emph{Comm. Pure Appl. Math.} \textbf{51} (1998), 23--45.

\bibitem{bib18}
Y.-H.~Kim and B.~Pass,
Multi-marginal optimal transport on Riemannian manifolds,
\emph{Amer. J. Math.} \textbf{137} (2015), 1045--1060.

\bibitem{bib16}
Y.-H.~Kim and B.~Pass,
Wasserstein barycenters over Riemannian manifolds,
\emph{Adv. Math.} \textbf{307} (2017), 640--683.

\bibitem{bib20}
S.~Kum, M.~H.~Duong, Y.~Lim and S.~Yun,
A GPM-based algorithm for solving regularized Wasserstein barycenter
problems in some spaces of probability measures,
\emph{J. Comput. Appl. Math.} \textbf{416} (2022), 114588.

\bibitem{bib8}
T.~Le~Gouic and J.-M.~Loubes,
Existence and consistency of Wasserstein barycenters,
\emph{Probab. Theory Related Fields} \textbf{168} (2017), 901--917.

\bibitem{bib2}
R.~J.~McCann,
A convexity principle for interacting gases,
\emph{Adv. Math.} \textbf{128} (1997), 153--179.

\bibitem{bib17}
B.~Pass,
On the local structure of optimal measures in the multi-marginal
optimal transportation problem,
\emph{Calc. Var. Partial Differential Equations}
\textbf{43} (2012), 529--536.

\bibitem{bib7}
B.~Pass,
Optimal transportation with infinitely many marginals,
\emph{J. Funct. Anal.} \textbf{264} (2013), 947--963.

\bibitem{bib22}
G.~Peyr\'e and M.~Cuturi,
Computational optimal transport,
\emph{Found. Trends Mach. Learn.} \textbf{11} (2019), 355--607.

\bibitem{bib9}
J.~Rabin, G.~Peyr\'e, J.~Delon and M.~Bernot,
Wasserstein barycenter and its application to texture mixing,
in \emph{Scale Space and Variational Methods
in Computer Vision},
Lecture Notes in Computer Science, vol.~6667,
Springer, Berlin, Heidelberg, 2012, pp.~435--446

\bibitem{ROF92}
L.~I.~Rudin, S.~Osher and E.~Fatemi,
Nonlinear total variation based noise removal algorithms,
\emph{Physica D} \textbf{60} (1992), 259--268.

\bibitem{San15}
F.~Santambrogio,
\emph{Optimal Transport for Applied Mathematicians:
Calculus of Variations, PDEs, and Modeling},
Progress in Nonlinear Differential Equations and Their Applications,
vol.~87, Birkh\"auser, Cham, 2015.



\bibitem{Sion58}
M.~Sion,
On general minimax theorems,
\emph{Pacific J. Math.} \textbf{8} (1958), 171--176.

\bibitem{bib21}
S.~Srivastava, C.~Li and D.~B.~Dunson,
Scalable Bayes via barycenter in Wasserstein space,
\emph{J. Mach. Learn. Res.} \textbf{19} (2018), 1--35.

\bibitem{bib3}
K.-T.~Sturm,
Probability measures on metric spaces of nonpositive curvature,
\emph{Contemp. Math.} \textbf{338} (2003), 357--390.

\bibitem{bib12}
C.~Villani,
\emph{Topics in Optimal Transportation},
Graduate Studies in Mathematics,
American Mathematical Society, Providence, RI, 2003.

\bibitem{bib11}
C.~Villani,
\emph{Optimal Transport: Old and New},
Grundlehren der mathematischen Wissenschaften,
vol.~338, Springer, Berlin, Heidelberg, 2009.


\end{thebibliography}
\end{document}